\documentclass[11pt]{article}

\usepackage{epsf}
\usepackage{psfrag}
\usepackage{amsmath}
\usepackage{amssymb}
\usepackage{lscape}
\usepackage{latexsym}
\usepackage[usenames]{color}
\usepackage[mathcal]{eucal}
\usepackage{stmaryrd}
\usepackage{textcomp}
\usepackage{setspace} 
\usepackage{graphicx} 
\usepackage{pgfplots} 
\usepackage{epsfig}   
\usepackage{epstopdf} 
\usepackage{appendix}
\usepackage{afterpage}
\usepackage{rotating} 
\usepackage{float} 
\usepackage[misc,geometry]{ifsym} 
\usepackage{comment}

\newtheorem{casestudy}{Case Study}[section]

\newcommand{\ds}{\displaystyle}
\newcommand{\dZ}{{\cal Z \kern -0.7em Z}}
\newcommand{\dC}{{\rm\hbox{C \kern-0.8em\raise0.2ex\hbox{\vrule height5.4pt width0.7pt}}}}
\newcommand{\dQ}{{\rm\hbox{Q \kern-0.85em\raise0.25ex\hbox{\vrule height5.4pt width0.7pt}}}}

\newcommand{\proofbox}{\hspace{\fill}{$\Box$}}

\newcommand\old[1]{}
\newcommand{\beqa}{\begin{eqnarray*}}
\newcommand{\eeqa}{\end{eqnarray*}}

\begin{document}

\title{\bf \Large{A Multiobjective Water Allocation Model for Economic Efficiency and Environmental Sustainability: Case Study}}
\author{Nahid Sultana \thanks{Department of Business Administration, International Islamic University Chittagong, {\tt email: nsultana@iiuc.ac.bd}}
\and
M. M. Rizvi \thanks{Centre for Smart Analytics (CSA), Institute of Innovation, Science and Sustainability, Federation University Australia, {\tt email: rizvimath@gmail.com }}
\and I.B. Wadhawan \thanks{The  Queensland University of Technology, {\tt email: 	indu.wadhawan@qut.edu.au}.(corresponding author)}}
\maketitle
\noindent \textbf{Abstract}  The management of irrigation water systems has become increasingly complex due to competing demands for agricultural production, groundwater sustainability, and environmental flow requirements, particularly under hydrologic variability and climate uncertainty. Addressing these challenges requires optimization frameworks that can jointly determine optimal crop allocation, groundwater pumping, and environmental flow releases while maintaining economic and hydrological feasibility. However, existing hydro-economic models, including the widely used Lewis and Randall formulation, may overestimate net benefits by allowing infeasible negative pumping and surface water allocations. We extend the Lewis and Randall framework by reformulating groundwater pumping and surface water use as non-negative, demand-driven decision variables and by explicitly incorporating environmental flow and canal capacity constraints. Three models are developed to maximize economic benefit, minimize environmental deficits, and a multiobjective model that evaluates the trade-offs between
these two objectives. An illustrative test case examining optimal crop area allocation and environmental flow management across dry, average, and wet years, using data from the Rajshahi Barind Tract in northwestern Bangladesh, is presented. The results show that the proposed formulation produces economically and hydrologically consistent solutions, identifying optimal strategies when either net benefits or environmental protection is prioritized, as well as Pareto-optimal trade-offs when both objectives are considered together. These findings provide practical insights for balancing farm income, groundwater sustainability, and ecological protection, offering a robust decision-support tool for irrigation management in water-limited river basins.


\textbf{Keywords} Irrigation Water Allocation, Environmental Flows, Crop Allocation, Multiobjective Optimization, Muhuri Irrigation Project.

\textbf{AMS} subject classifications. 90B90, 90C29, 90C30, 92B05

\section{Introduction}\label{sec1}

\noindent Water is the foundation of agricultural production, yet its availability is becoming increasingly uncertain due to climate variability, population growth, and rising competition across sectors. Irrigation remains the largest consumer of freshwater globally, accounting for nearly 70\% of total withdrawals \cite{Taye2021}. In developing regions dominated by monsoon climates, water availability fluctuates sharply between wet and dry seasons, complicating equitable and sustainable allocation. Balancing economic efficiency with ecological integrity has therefore emerged as one of the most critical challenges in contemporary water-resources planning and modelling \cite{Jeuland2013}.

\noindent Across South Asia, river basins are experiencing intensified competition among agriculture, domestic demand, and environmental needs. Bangladesh represents a particularly vulnerable case, where extreme monsoon inflows are followed by severe dry-season scarcity \cite{Ullah2021}. In the northwestern region, especially the Rajshahi Barind Tract, irrigated agriculture faces acute stress due to declining river flows, erratic rainfall, and widespread groundwater depletion. Farmers increasingly rely on shallow and deep tube wells to supplement canal supplies, which has improved short-term productivity but raised serious long-term sustainability concerns \cite{Dey2017}. Effective water allocation in this context requires modelling frameworks that explicitly capture hydrological variability, competing sectoral objectives, and environmental flow requirements under local socio-economic and institutional constraints.

\noindent Conventional allocation policies often prioritize crop yield or economic returns while insufficiently accounting for environmental flow requirements essential for sustaining riverine ecosystems and groundwater recharge \cite{Tulip2022}. Conversely, strictly conservation-driven approaches can constrain agricultural output and threaten rural livelihoods. This inherent trade-off highlights the need for optimization-based decision frameworks capable of simultaneously addressing economic and ecological objectives \cite{Xevi2005}. Multiobjective optimization, in particular, provides a systematic means to explore such trade-offs and identify Pareto-efficient solutions that reflect realistic policy alternatives.

\noindent Environmental flow (e-flow) is defined as the quantity, timing, and quality of water required to maintain riverine ecosystems and the human livelihoods that depend on them (Brisbane Declaration, 2007). Rather than being a fixed threshold, environmental flows vary seasonally in response to natural hydrological regimes \cite{Pastor2014}. These variations sustain ecological processes such as fish migration, floodplain renewal, and wetland regeneration. Excessive abstraction for irrigation disrupts this natural rhythm, resulting in ecosystem degradation, reduced groundwater recharge, and long-term socio-economic consequences. In planning terms, deviations below a defined target environmental flow are quantified as Environmental Flow Deficiency (EFD), which serves as a direct indicator of ecological stress.

\noindent Target Environmental Flow (TEF) represents the balance point between human water use and ecosystem needs and is commonly expressed as a proportion of monthly inflow or through seasonal hydrographs \cite{Poff2010}. The choice of TEF formulation critically influences model outcomes. Fixed targets often exaggerate deficiencies, whereas seasonally adaptive targets allow more realistic system responses.
Numerous methods have been proposed to estimate environmental flow requirements, including holistic approaches such as the Building Block Methodology (BBM), variability-based methods like the Range of Variability Approach (RVA) and Variable Monthly Flow  (VMF), and hydrological rules such as Tennant and Smakhtin.
\par
\noindent Previous comparative studies evaluating the Smakhtin, Tennant, Tessmann, Variable Monthly Flow (VMF), and Q90-Q50 methods have shown that environmental flow estimates derived from the Tessmann and VMF methods exhibit the highest agreement with locally assessed environmental flow requirements (R² = 0.91) \cite{Pastor2014} reflecting their ability to capture intra-annual hydrological variability. Although the Tennant method remains one of the earliest and most widely applied approaches due to its simplicity and low data requirements, it relies on fixed fractions of mean annual flow and does not explicitly represent seasonal dynamics. In this study, the Tessmann method is adopted because it provides stricter protection during low-flow periods by allocating a higher environmental share, while allowing proportional relaxation during high-flow months. This precautionary yet seasonally adaptive structure is well suited to the pronounced dry-season scarcity, groundwater dependence, and ecological sensitivity of the Rajshahi Barind Tract.

\noindent The integration of economic and ecological objectives within irrigation allocation models has evolved steadily over recent decades. Early work by Xevi and Khan \cite{Xevi2005} demonstrated that multiobjective optimization could reconcile agricultural productivity with ecosystem preservation under semi-arid conditions. Subsequent studies incorporated environmental flows directly into allocation frameworks \cite{Meijer2012}, applied evolutionary algorithms under uncertainty \cite{Sadati2014,Li2014}, and explored Pareto-optimal trade-offs between income and hydrological stability \cite{Lalehzari2016,Mohammadrezapour2019}. Applications across Australia \cite{Adamson2005, Lewis2017}, South Africa \cite{Ikudayisi2018}, Europe \cite{Crespo2019}, and Iran \cite{Zeinali2020} further highlighted the importance of explicitly modelling trade-offs under variable climatic and institutional settings.

\noindent Recent research has increasingly emphasized sustainability-oriented perspectives. Studies by Marzban et al. \cite{Marzban2021} and Bhuiyan \cite{Bhuiyan2022} underscored that many allocation models still underestimate ecological costs, leading to persistent imbalances between economic objectives and environmental safeguards. Ullah and Micah \cite{Ullah2021} demonstrated that rainfall variability plays a dominant role in shaping both net benefits and environmental flow deficiency in Bangladesh, reinforcing the need to jointly consider rainfall–inflow interactions. However, recent studies focused on the Rajshahi region have largely relied on direct, linear, and deterministic analyses of cropping patterns, rainfall trends, and irrigation demand, often neglecting resource constraints, environmental flow deficiency, and optimized water-allocation strategies \cite{ HAli2021,Hossain2020, Masum2013,  Mila2016,Milla2018, Mojid2020, Rahman2015, Rashid2017,Sen2019},

\noindent Comparative assessments of environmental-flow estimation methods reveal substantial variation in recommended allocations. While Tennant’s approach simplifies calculations using mean annual flow \cite{Tennant1976}, methods like Tessmann allocate higher environmental shares during low-flow months \cite{Tessmann1980}, potentially constraining irrigation in water-scarce systems. Variability-based approaches, such as Variable Monthly Flow (VMF) adapt allocations to seasonal hydrological conditions and are often better suited to monsoon-dominated basins. Integrating such dynamic environmental requirements into optimization-based irrigation models remains an open research challenge, particularly in data-limited regions.

\noindent Despite growing recognition of hydrological variability, many allocation studies still rely on average-year or deterministic assumptions \cite{Poff2010}, overlooking interannual variability and conjunctive surface–groundwater use \cite{Siebert2010}. In regions like Rajshahi—where declining surface flows intensify groundwater abstraction—such omissions risk misrepresenting sustainability thresholds \cite{Islam2020,Tabassum2025}. Although South Asian studies have examined large transboundary systems such as the Ganges and Brahmaputra, regional irrigation subsystems remain underexplored.

\noindent The Rajshahi Barind Tract exemplifies the urgent need for multiobjective, scenario-based irrigation allocation models that integrate hydrological uncertainty, environmental flow requirements, and economic decision-making under local institutional constraints. While Sultana et al. \cite{Sultana2025} advanced this direction by jointly maximizing net benefit and minimizing environmental flow deficiency, they highlighted the need for further refinement through crop-specific land constraints, canal-capacity limitations, and seasonal environmental-flow targets. Addressing these gaps is essential for developing decision-support tools capable of guiding sustainable irrigation management in water-stressed regions such as Rajshahi.

\par
\par
\par
\subsection{Research Gap and the New Contribution}
Many existing models used to plan crop production and irrigation water use, such as the model by Lewis and Randall (2017), calculate net profits based on crop area, crop water needs, and groundwater pumping. However, these models allow groundwater pumping and surface water use to become negative when crops do not need water. This creates unrealistic results and often leads to an overestimation of net income. For example, when crops require no water but surface water is still available, the model produces negative pumping values and incorrectly treats them in the objective function as income rather than a cost. A similar problem occurs when more surface water is available than crops actually need, resulting in negative surface-water costs.

\noindent In addition, many models do not include physical limits of the water system, such as the capacity of canals to carry water. This allows water allocations higher than the irrigation system can actually deliver.

\noindent Data limitations also create a major challenge, particularly in the Rajshahi Barind Tract, where key information on crops, water use, and inflows is often missing or incomplete and must be carefully collected and processed. At the same time, solving these large, complex models requires efficient algorithms. Because the models include many constraints and rely heavily on data, not all solvers can handle them effectively, so it is necessary to test and compare different solution methods.

\noindent Therefore, there is a clear need for a new modeling framework that ensures groundwater pumping and surface water use are always non-negative and based on real crop water demand, that respects canal capacity and environmental flow limits, and that is supported by efficient and well-tested solution algorithms. Filling this gap is essential to produce realistic estimates of farm profits and to design water management policies that are both economically sensible and physically achievable.



\noindent The remainder of the paper is organized as follows. Section \ref{Insight} provides actionable decision and managerial insights. Section \ref{Case Study in Rajshahi} presents the case description along with illustrations of the relevant crop and weather data for the Rajshahi region. Section \ref{Moo} introduces the key concepts underlying the three developed models. Section \ref{NuExp} reports the numerical experiments, presents the results, and discusses the major findings, with Section \ref{MOP} detailing the resulting Pareto front and comparing solutions across different hydrological conditions. Finally, Section \ref{Con} summarizes the key findings, discusses limitations, and outlines directions for future research.

\section{Decision and Managerial Insights} \label{Insight}

\noindent 
Our study presents a structured decision-making framework for agricultural water allocation that integrates both economic and environmental considerations under realistic hydrological conditions. The proposed model addresses several key challenges often faced in irrigation management with limited canal capacity, crop-specific minimum land requirements, and the joint use of surface and groundwater resources. By introducing these constraints into a multiobjective optimization framework, the model allows policymakers to evaluate trade-offs between maximizing crop benefits and maintaining ecological flow requirements across dry, average, and wet years.
\par
\noindent The model consists of two interrelated components: one dealing with economic allocation efficiency and the other with environmental sustainability. The first component focuses on optimizing water distribution among crops to achieve the highest total net benefit while respecting canal capacity and groundwater pumping limits. This helps irrigation planners identify the most effective combination of crops and water sources, ensuring that resources are neither wasted nor overexploited. The second component incorporates environmental flow considerations, balancing water withdrawals with ecosystem needs and maintaining the health of downstream habitats.
\par
\noindent From a managerial perspective, the framework provides valuable insights for water authorities, policymakers, and farm cooperatives. It enables them to make informed decisions about how much surface water to allocate, how much groundwater to pump, and how to adjust these allocations when inflow conditions change. The model’s use of Sequential Quadratic Programming (SQP) ensures fast and accurate solutions, making it suitable for real-time or seasonal planning scenarios. In practical application, the framework can guide canal scheduling, optimize seasonal cropping patterns, and support equitable distribution strategies under limited resource availability.
\par
\noindent Beyond its immediate operational benefits, the model also contributes to long-term sustainability goals. Explicitly linking environmental flow requirements with economic objectives, it helps decision-makers understand the consequences of overuse or misallocation, not only in terms of reduced yields but also in ecological degradation. Sensitivity analysis further clarifies how variations in canal capacity, inflow volume, or crop area requirements affect both total benefit and environmental compliance. These insights support adaptive planning and strengthen the resilience of irrigation systems against climatic variability.
\par
\noindent In essence, this research bridges the gap between theoretical optimization and practical irrigation management. It equips planners and decision-makers with a transparent and flexible tool to allocate water resources efficiently, sustain agricultural productivity, and protect the environment, a crucial step toward achieving integrated water resource management in water-stressed regions like Rajshahi.

\section{Case Description: Rajshahi Barind Tract } \label{Case Study in Rajshahi}
\par
\noindent The Rajshahi Barind Tract in northwestern Bangladesh is selected as the case study due to its pronounced hydrological variability, strong dependence on irrigation, and the increasing conflict between agricultural water demand and environmental flow requirements. Agricultural production in the region follows a distinct seasonal pattern characterized by severe water scarcity during the low-flow months and relatively abundant water availability during the monsoon season.  In this study, a total of nine major crops are considered, of which six are predominantly cultivated during the low-flow period, when irrigation demand is highest and competition between agricultural use and environmental flow protection is most critical. From May to October, rainfall and river inflow increase substantially, and farmers in the Rajshahi region typically cultivate very short-duration local crops, for instance, vegetables, Mung bean, Black gram, sesame, etc. Because these crops have brief growing periods, flexible planting schedules, and rapidly changing production structures, they are not explicitly included in the optimization framework. Monthly rainfall, reference evapotranspiration, and river inflow exhibit strong intra-annual variability, resulting in surplus surface water during the monsoon and acute shortages during the dry season. Target environmental flows are defined as a proportion of monthly inflow to preserve riverine ecological functions, further constraining surface water availability during low-flow periods. Previous studies in the Rajshahi region have primarily relied on direct, deterministic analyses of cropping patterns, rainfall trends, and irrigation demand, often overlooking key resource constraints, environmental flow deficiency, and optimized water allocation decisions. In contrast, the present study incorporates multiple crops with distinct water requirements, crop calendars, yields, and heterogeneous economic characteristics within a multiobjective optimization framework. While the total cropland area is fixed, the effective cultivable land allocated to each crop varies with productivity per hectare. Groundwater abstraction is limited by an upper pumping constraint, and surface water diversion is restricted by canal capacity. Together, these hydrological, agronomic, and infrastructural constraints create a complex allocation problem that requires balancing economic returns with environmental flow protection, particularly in dry and average hydrological years.

\subsection{Calendar of selected crops by growth Stage}\label{crop calendar} 
\begin{figure}[H]
\hspace{-1cm}
\begin{center}
\includegraphics[width=180mm]{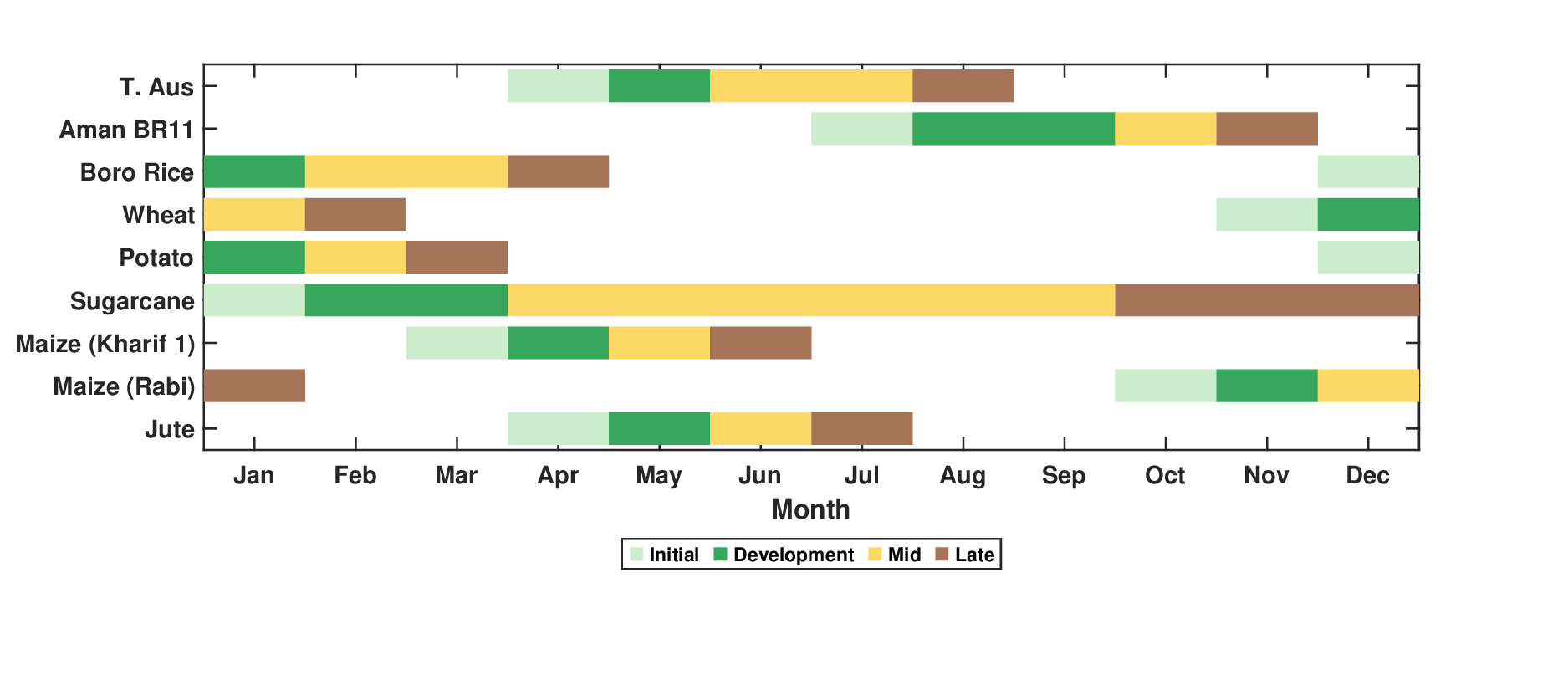} 
\end{center}
\caption{Crop Calendar by Growth Stage in Rajshahi.
}
\label{Crop Calendar}
\end{figure}

\noindent Figure~\ref{Crop Calendar} illustrates the cropping calendar for the major profitable crops in Rajshahi, northwest Bangladesh, showing the timing of each crop across four growth stages: initial, development, mid, and late. Each stage has a different water requirement, expressed through crop coefficients that were used to estimate irrigation demand in the model. The seasonal rice varieties (T. Aus or Aus, Aman BR11, and Boro) spread across different months from March to December, while wheat and potato grow mainly in the dry season between November and April. Sugarcane lasts almost the whole year because of its long growth cycle. Jute and maize (both Kharif 1 and Rabi types) reflect the region’s typical monsoon and winter crops.
\subsection{Crop Production}\label{in Rajshahi (Tons/ha)} 
\begin{figure}[H]
\hspace{-1cm}
\begin{center}
\includegraphics[width=120mm]{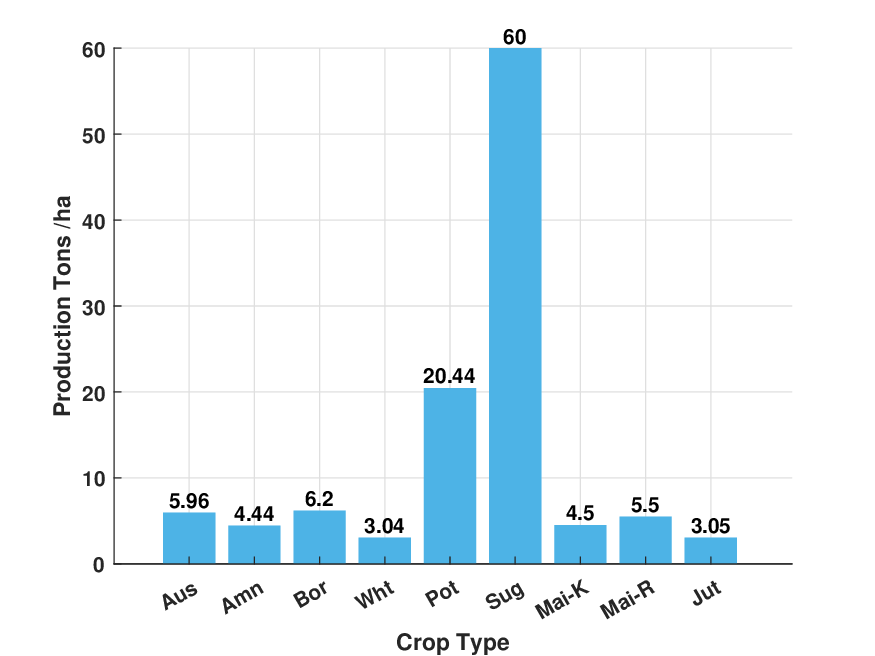} \\
\hspace*{-2.9cm}
\caption{Crop Production in Rajshahi (Tons/ha)}
\label{Crop Production}
\end{center}
\end{figure}
\noindent Figure~\ref{Crop Production} depicts the variation in crop yields across different types grown in the Rajshahi region, measured in tons per hectare.
Overall, the pattern highlights that Rajshahi’s crop production is heavily tilted toward high-yield and water-intensive crops like sugarcane and potato, whereas cereals such as wheat and jute contribute modestly to the total agricultural output.

\subsection{Inflow and Target Environmental Flow(TEF)}\label{Inflow and Target Environmental Flow(TEF)} 
\begin{figure}[H]
\hspace{-1cm}
\begin{center}
\includegraphics[width=150mm]{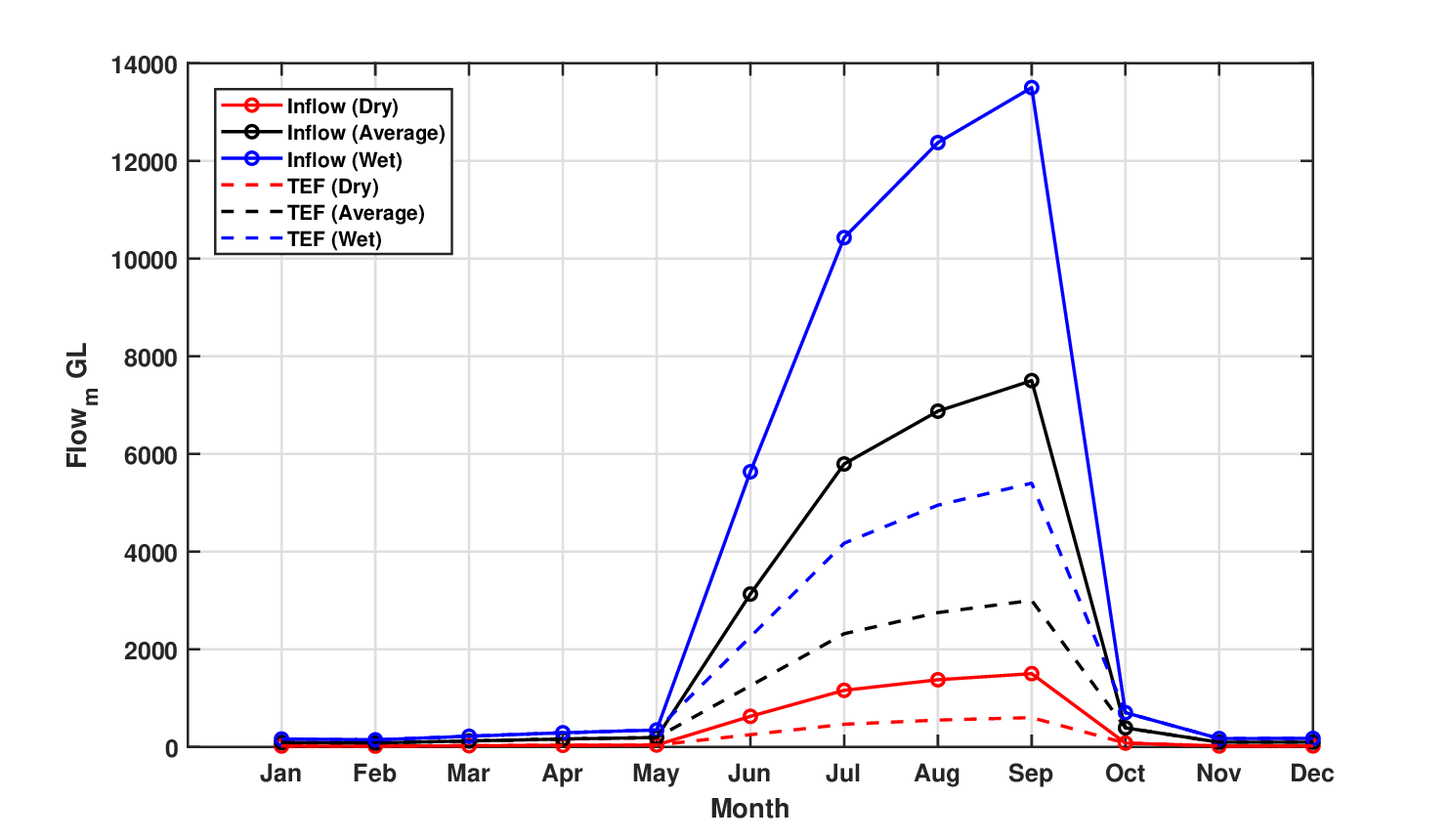} \\
\end{center}\hspace{-1cm}
\vspace{-1cm}
\caption{ Inflow and TEF in Dry, Avg. and Wet years}
\label{Inflow and TEF in Dry, Avg. and Wet years}
\end{figure}

\noindent Figure ~\ref{Inflow and TEF in Dry, Avg. and Wet years} shows the monthly variation of inflow and target environmental flow (TEF) under dry, average, and wet hydrological conditions. Following the Tessmann method, TEF equals 100\% of inflow during low-flow months (November–April), ensuring full ecological protection. From May onward, inflow rises sharply and peaks around September, especially in wet years. During these high-flow months, TEF is defined as a proportion of monthly inflow, allowing up to 60\% of surface water to be allocated for irrigation. This approach protects ecological needs in dry periods while providing greater allocation flexibility during high-flow months.

\subsection{Rainfall And Evapotranspiration}
\begin{figure}[H]
\hspace{-1cm}
\begin{center}
\includegraphics[width=150mm]{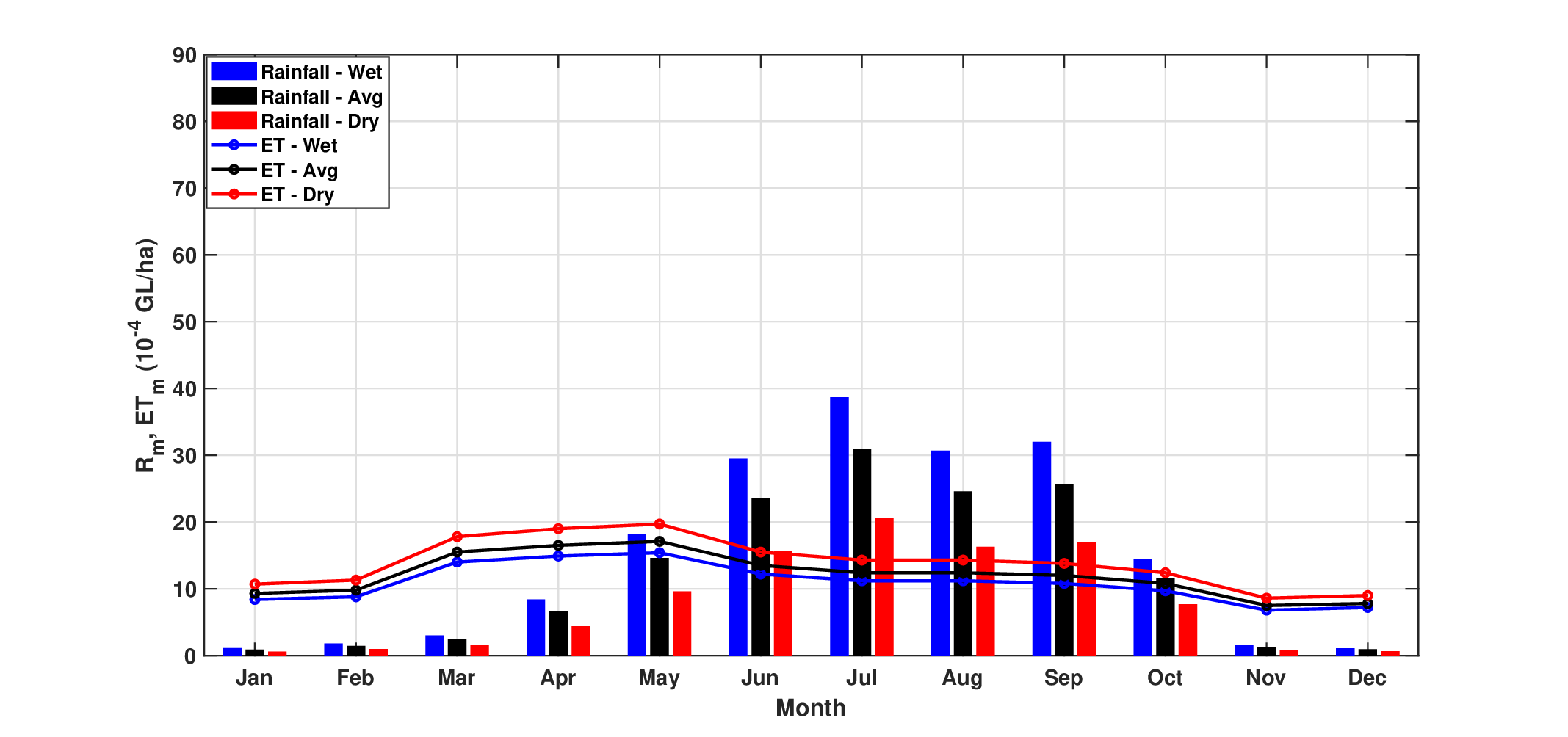} \\
\caption{Rainfall($R_m$) and  Evapotranspiration($ET_m$) in Dry, Avg. and Wet years.}
\label{fig_R_m & ET_m}
\end{center}\hspace{-1cm}
\end{figure}
The combined rainfall and evapotranspiration pattern in Figure~\ref{fig_R_m & ET_m}  illustrates the seasonal water balance across dry, average, and wet years in the study area. Rainfall, shown as bars, remains minimal from January to March, increases sharply from April, and peaks during the monsoon months of June to August, before declining rapidly after September. Wet years record markedly higher rainfall throughout the season, whereas dry years show a pronounced deficit during the same period. In contrast, the reference evapotranspiration (ET), represented by lines, rises gradually from February, reaching its maximum during April–May when temperature and solar radiation are highest, and then decreases during the monsoon due to increased humidity and cloud cover. ET remains lowest during wet years and highest during dry years, indicating the influence of climatic variation on atmospheric water demand. Overall, the figure highlights the inverse relationship between rainfall and evapotranspiration and provides a basis for assessing seasonal irrigation requirements and hydrological performance under different climatic conditions.
\newpage
\section{Formulation the Mathematical Model and Solution Approach}\label{Moo}

\noindent 
Decision variables:
\begin{table} [H]
\begin{tabular}{l l}
$X_c$& Crop Area for each crop c in hectare $ha$\\
$Env.Flow_m$& Environmental Flow in month m in gigaliter $GL$\\
\end{tabular}
\end{table}
\noindent Input parameters:
\begin{table} [H]
\begin{tabular}{l l}
$P_c$& Price of total crop production (Tk./Ton)\\
$Y_c$& Crop production per Hectare (Ton/ ha)\\
$C_w$& Per unit Cost of Surface Water conveyance (Tk./GL)\\
$K_{cm}$  & Different Crop coefficients in the several growth stages of selected crops c in month m\\
$ET_m$& Monthly Reference Evapotranspiration($GL/ha$) \\
$R_m$  & Monthly Rainfall($GL/ha$) \\
$Allocation_m$& Surface Water allocated ($GL$) month wise\\
$P_m$  & Pumping Ground Water($GL$) month-wise\\
$C_p$& Per unit Cost of pumping Groundwater(Tk./GL)\\
$Vcost_c$& all related Variable costs for each crop c per area (Tk./ha)\\
$Inflow_m$& Available River water ($GL$) month wise \\
$Tar.Env.Flow_m$ & The minimum required Environmental flow($GL$) in each month \\
$T_{\text{pump}}$ & Total Pumping Water  $(GL)$\\
$T_{\text{area}}$ & Total Area $(ha)$\\
$min_{area}$ & Minimum Area $(ha)$
\end{tabular}
\end{table}

\noindent Decision variables dependent parameters:
\begin{table} [H]
\begin{center}
\begin{tabular}{l l}
$Allocation_m= max({Inflow}_m -  {Env.Flow}_m,0)$\\
$ P_m = max (\sum_{c} (K_{cm} \cdot ET_m - R_m  ) X_c - Allocation_m$,0)\\
\end{tabular}
\end{center}
\end{table}
\noindent The objective functions developed by Lewis and Randall\cite{Lewis2017} is as follows:
\par
\noindent The first objective, $f_1$:
\begin{multline}\label{revobj}
    NB_{\text{max}} =  \sum_{c} P_{c}Y_c X_c - C_w \cdot  \sum_{m} \left(\sum_{c} (K_{cm} \cdot ET_m -R_m) X_c -  P_m \right)\\
    - C_p \cdot \sum_{m} P_m - \sum_{c} Vcost_c \cdot X_c
\end{multline}

\noindent The economic objective in \eqref{revobj} can result in overstated net benefits under certain conditions. Specifically, if in any month the crop water requirement is zero that is, \( \sum_{c} (K_{cm} \cdot ET_m -R_m) X_c  = 0 \) for some $m$, but in the same month some surface water is available for irrigation implies that the allocation is positive (\(Allocation_m > 0\)). Therefore, the pumping definition produces a negative value (\( P_m < 0 \))  and thus is incorrectly counted as revenue in \eqref{revobj} rather than being treated as a cost in irrigation management. This situation is illustrated in Figure 7, where from June to October the crop water requirement is zero, but surface water available in these months, this trigger the negative pumping water. To address this issue, the pumping term is constrained in \eqref{revobj} using \( P_m = \max \left\{ 0, \sum_{c}(K_{cm} ET_m - R_m )X_c - \text{Allocation}_m \right\} \), ensuring that pumping occurs only when crop water demand exceeds surface water allocation. This modification enables groundwater pumping costs to reflect actual irrigation needs and provides a more accurate estimate of net benefits.

\noindent A similar overestimation of revenue objective in \eqref{revobj} can occur when calculating surface water costs. To address this, surface water allocation is constrained to be non-negative, ensuring that available inflow is greater than or equal to the required environmental flow that is, $Inflow_m \ge Env\cdot Flow_m$ for all $m$. Incorporating this condition directly into the model constraints restricts allocations to non-negative values, resulting in a more accurate estimation of surface water costs.
\noindent The extended objective function can thus be written as:
\begin{multline}\label{revobj1}
    f_1(X_c, Env.Flow_m) =  \sum_{c} P_{c}Y_c X_c - C_w \cdot  \sum_{m} \left(\sum_{c} (K_{cm} \cdot ET_m -R_m) X_c -  P_m \right)\\
    - C_p \cdot \sum_{m} \max \left\{ 0, \sum_{c}(K_{cm} ET_m - R_m )X_c - \text{Allocation}_m \right\} - \sum_{c} Vcost_c \cdot X_c
\end{multline}

\noindent In the proposed model, we also incorporate canal constraints, defined as
\[
Inflow_m - Env.Flow_m \le C_c, \;\; \forall m.
\]
to ensure that water flow remains within the feasible capacity of the canal.

\noindent Hence, the proposed model for maximizing revenue subject to constraints is as follows:

\noindent {\bf Model 1: Optimizing the Net Benefit (NB)}
\begin{equation}\label{eqmathemodel1}
\begin{array}{rl} 
\ds \max_{(X_c, Env.Flow_m) \in X} & \ f_1(X_c, Env.Flow_m)\\
\mbox{subject to the constraints} \\
& \ds X_c, Env.Flow_m \geq 0.
\end{array}
\end{equation}
where, \begin{equation*}
\hspace*{-5mm}
X := 
\left\{
\ds (X_c, \text{Env.Flow}_m) \in \mathbb{R}^{c} \times \mathbb{R}^{m} \;\middle|\;
\begin{array}{l}
\ds \sum_{m} \max \Big\{ 0, \sum_{c}(K_{cm} ET_m - R_m )X_c - \text{Allocation}_m \Big\} \leq T_{\text{pump}}, \\[1ex]
\ds \text{min}_{\text{area}} \leq X_c, \quad \sum_{c} X_c \leq T_{\text{area}}, \\[0.5ex]
\text{Env.Flow}_m  \leq \text{Inflow}_m, \;\;\forall m, \; \text{Inflow}_m - \text{Env.Flow}_m \le C_c, \;\; \forall m.
\end{array}
\right\}.
\end{equation*}

The second objective, $f_2$:
\begin{equation}\label{envobj}
    EFD_{\text{min}} = \sum_{m} \max \left( \text{Tar.Env.Flow}_m - \text{Env.Flow}_m,\ 0 \right)
\end{equation}
\noindent The model includes two decision variables, $X_c$ and $Env.Flow_m$. The objective is to maximize net benefit while limiting adverse environmental impacts. The model is parameterized by $P_c$, $Y_c$, $C_w$, $K_{c,m}$, $ET_m$, $R_m$, $T_{\text{pump}}$, $Allocation_m$, $T_{\text{area}}$, $min_{\text{area}}$, $P_m$, $C_p$, $Inflow_m$, $C_c$ and $Vcost_c$. The resulting trade-off analysis reveals that the corner solution, which achieves the maximum economic benefit, corresponds to the weakest environmental performance, clearly illustrating the tension between profit maximization and environmental protection.
\par
\noindent {\bf Model 2: Optimizing the EFD}
\begin{equation}\label{eqmathemodel2}
\begin{array}{rl} 
\ds \min_{(X_c, Env.Flow_m) \in X} & \ f_2(X_c, Env.Flow_m)\\
\mbox{subject to the constraints} \\
& \ds X_c, Env.Flow_m  \geq 0.
\end{array}
\end{equation}
\noindent Accordingly, the model places primary emphasis on minimizing Environmental Flow Deficiency (EFD). When the prescribed target environmental flow does not exceed the realized environmental flow, the EFD reduces to zero. This strict constraint ensures full compliance with environmental flow requirements; however, it does so by sacrificing economic performance, resulting in the lowest achievable level of net benefit in favor of strong environmental protection.
\par
\noindent {\bf Model 3: Trade-off between NB and EFD }
\begin{equation}\label{eqmathemodel3}
\begin{array}{rl} 
\ds \min_{(X_c, Env.Flow_m) \in X}  & \ds \ \left\{-f_1, f_2\right\}\\
\mbox{subject to the constraints} \\
& \ds X_c, Env.Flow_m \geq 0.
\end{array}
\end{equation}
\noindent The proposed multiobjective programming model simultaneously considers two conflicting objectives: minimizing the negative net benefit function, $f_1$ and minimizing the environmental flow deficiency function, $f_2$. The constraint framework incorporates feasible pumping limits and minimum land-use requirements, while no upper bound is imposed on crop land allocation in order to allow full flexibility in land-use decisions. In addition, non-negativity constraints ensure that crop areas and environmental flow allocations remain physically meaningful. Collectively, this formulation explicitly captures the inherent trade-off between maximizing agricultural economic returns and preserving essential ecological flow requirements. By representing the complete set of non-dominated solutions, the model is able to provide a comprehensive view of the Pareto-optimal trade-off space.

\noindent The models are implemented in MATLAB using its standard solvers. The step-by-step solution approach and the algorithms applied to Models 1–3 are provided in Appendices \ref{solprocedure} and \ref{mopalg}.

\noindent The following assumption is used to ensure the ease of the developed model.
\begin{itemize}

\item Since all the lands fall within the same Barind region, it is assumed that the climatic and edaphic conditions, namely temperature, evapotranspiration, humidity, rainfall, soil moisture, and soil fertility are homogeneous across the study area within each representative year type (dry, average, and wet).
\end{itemize}
\section{Numerical Experiments,  Results and Sensitivity Analysis}\label{NuExp}

\noindent A numerical experiment is conducted for three representative hydrological conditions, such as dry, average, and wet years, to optimize irrigation water allocation under multiple resource constraints arising from environmental, social, and governance considerations. Two single-objective formulations are first analyzed in Case Study \ref{exp1}. In Model~1 in ~\eqref{eqmathemodel1}, the net benefit is maximized, and Model~2 in \eqref{eqmathemodel2}, which minimizes environmental flow deficiency. In Case Study \ref{exp2}, we present the numerical implementation of Model~3 in \eqref{eqmathemodel3}, which approximates the non-dominated front representing the trade-off between maximizing profit and minimizing environmental flow deficiency. The optimization problems are implemented in MATLAB using the `fmincon' solver. The interior-point, active set, and sequential quadratic programming (SQP) algorithms are employed to ensure robustness under a large number of constraints and to assess solver performance for large-scale nonlinear optimization. Computational efficiency is evaluated based on how well optimal solutions are approximated, with performance depending on the algorithms selected in fmincon and the initial guesses provided. For each solver, the computational time and the quality of the obtained solutions are reported. All numerical experiments are carried out on an HP Evo laptop equipped with 16 GB of RAM and a Core i7 processor operating at 4.6 GHz.
\par
\begin{casestudy}\label{exp1}
This numerical example is based on the case study in the Rajshahi zone. The study considers a regional water board managing nine crops, including staple foods and other highly demanded crops, over 36 months across three different hydrological years. The total cultivable area is 182,271 hectares, with a maximum annual water pumping capacity of 500 GL and different minimum cropland as its demand. The costs for pumped water and surface water are 100,000 and 26,000 $Tk.$ $per$ $GL$, respectively. The target environmental flow is set following the Tessman method \cite{Tessmann1980}. For the low-flow months, all water should be reserved for environmental flows, while in high-flow months, this could be up to 40 percent of the inflow. Surface water is conveyed through the canals with a maximum carrying capacity of 6000 GL. Therefore, there is no fixed target environmental flow over the periods. Additional parameter data, including rainfall, evapotranspiration, crop coefficients, crop production, total crop income, and variable costs, are obtained from various sources, particularly academic journals and from the Tables ~\ref{table:ex50},\ref{table:ex5b},\ref{table:ex5c} in Appendix~\ref{appen}.
We examine how crops can be allocated subject to the minimum area requirement, with no upper bound of the variables imposed, so that the allocation remains fully flexible. To satisfy the minimum demand for each crop, we establish minimum cultivation areas. The required minimum areas for each crop are determined through surveys and incorporated into the model, as shown in Table \ref{table:mincrops}. The analysis considers whether the environmental flow meets its target requirement and evaluates cases with all the included five constraints such as, total pump limit, minimum area, total area, inflow and environmental flow, and canal capacity constrains in order to identify the maximum achievable profit. The results yield the individual optimal points of the objective functions, which are subsequently used in Case Study 2 to examine the trade-off between maximizing profit and minimizing environmental flow deficiency, thereby approximating the complete shape of the Pareto front. 
\end{casestudy}
\begin{table} [H]
 		\caption{\small{\textit{Minimum Allocation of Crops Per Year}.
        }}
		\footnotesize
		\vskip 1.5em
 		\centering
 		\begin{tabular}{|c| c| c| c| c| c| c| c| c| c| c}
 			\hline
 			crops &   Aus &  Aman  & Boro & Wheat& Potato & Suger & Maize & Maize & Jute   \\ 
             &   rice &  rice  & rice & &  & cane & Kharif 1& Rabi &    \\ [0.5ex]
 			\hline
$\text{min}_{\text{area}}$&20000 &35000 & 30000 &10000&15000&16000&5000&5000&6000\\ [.5ex]\hline

                \end{tabular}
		\label{table:mincrops}
	\end{table}
\begin{table} [H]
 		\caption{\small{\textit{Optimal Allocation of Crops and Environmental Flow with Pumped water for Model 1 in Dry Year. Optimal values are $f_1=-2.6746 \times 10^{10}$ and $f_2=194.9720$. }}}
		\footnotesize
		\vskip 1.5em
 		\centering
 		\begin{tabular}{|c| c| c| c| c| c| c| c| c| c| c| c|c|}
 			\hline
 			Crops &   Aus &  Aman  & Boro & Wheat& Potato & Suger & Maize & Maize & Jute &  & &  \\ 
             &   rice &  rice  & rice & &  & cane & kharif 1& rabi &  &  &  &  \\ [0.5ex]
 			\hline
$X_c$&20000 &35000 & 30000 &10000&55271&16000&5000&5000&6000& & &\\ [.5ex]\hline
Env.&Jan &Feb & Mar &April&May&June&July&Aug&Sep&Oct&Nov & Dec\\ [.5ex]\hline
Flow&0 &0 & 0 &0&0 &612&1159&1372&1500&50&0&0\\ [.5ex]\hline
$P_m$&69.85 &110.57 & 135.4 &74.79&26.59 &0&0&0&0&0&19.38&38.35\\ [.5ex]\hline
                \end{tabular}
		\label{table:ex1a}
	\end{table}
\begin{figure}[H]
\begin{center}
\includegraphics[width=140mm]{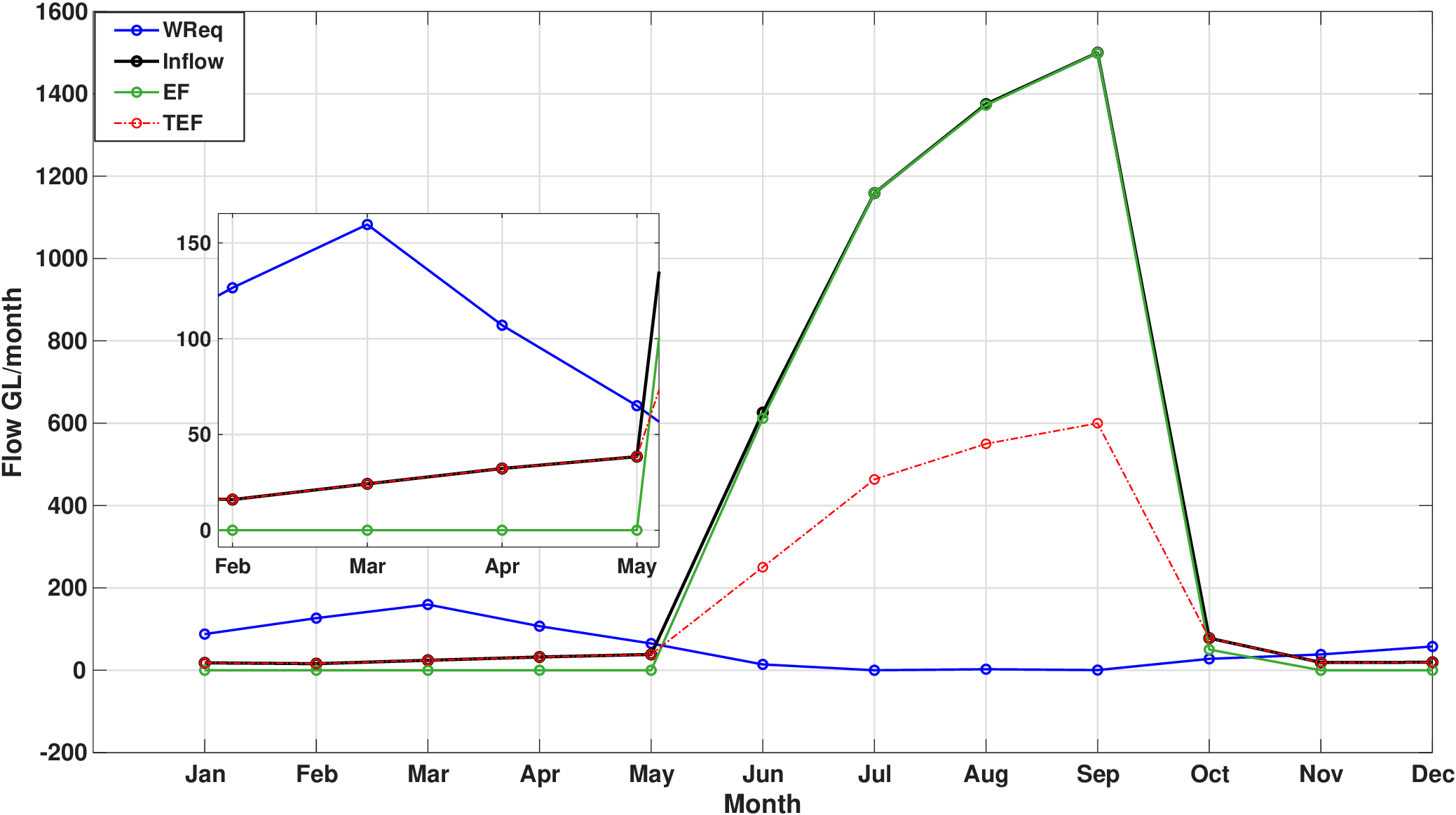} \\
\hspace*{-2.9cm}
{\scriptsize \caption{{Monthly Allocation of Environmental Flow and Water Requirement for Model 1 in Dry Year.
}}}
\end{center}
\label{fig_dry1}
\end{figure}
\noindent Table~\ref{table:ex1a} and Figure~\ref{fig_dry1}  illustrate the optimal crop and environmental flow allocation for Model~1 under the dry-season scenario, where the objective is to maximize net benefit. The model allocated the largest cultivation area to potatoes in order to maximize net benefit. This result is justified by the high production rates per hectare observed for both sugarcane and potato (see Figure \ref{Crop Production}). However, potatoes generate a total crop income of $13,000$ Tk per ton, which is more than double that of sugarcane ($5,580$ Tk per ton; see Table~\ref{table:ex50}). As a result, the model allocated $55,271$ hectares to potato cultivation to optimize overall income. We also observe that no environmental flow is allocated from November to May, as most crops are cultivated during this period and surface water availability is insufficient, requiring reliance on pumped water. The allocation of pumped water over twelve 12-month periods is presented in Table~\ref{table:ex1a}. Environmental flow releases increase sharply during the monsoon season and peak in September due to heavy rainfall. Under this allocation, the model achieves a maximum net benefit $2.6746 \times 10^{10}$ with a corresponding flow deficiency of $194.9720$ $GL$. 
\begin{table} [H]
 		\caption{\small{\textit{Optimal Allocation of Crops and Environmental Flow with Pumped water for Model 2  in Dry Year. Optimal values are
        $f_1=-1.7165 \times 10^{10}$ and $f_2=39.8460$.}}}
		\footnotesize
		\vskip 1.5em
 		\centering
 		\begin{tabular}{|c| c| c| c| c| c| c| c| c| c| c| c|c|}
 			\hline
 			crops &   Aus &  Aman  & Boro & Wheat& Potato & Suger & Maize & Maize & Jute &  & &  \\ 
             &   rice &  rice  & rice & &  & cane & kharif 1& rabi &  &  &  &  \\ [0.5ex]
 			\hline
$X_c$&20000 &35000 & 30000 &10000&15000&16000&5000&5000&6000& & &\\ [.5ex]\hline
Env.&Jan &Feb & Mar &April&May&June&July&Aug&Sep&Oct&Nov & Dec\\ [.5ex]\hline
Flow&18 &15 & 24 &27&9 &611&627&635&630&75&18&19\\ [.5ex]\hline
$P_m$&63.55 &77.16 & 108.5 &101.79&35.39 &0&0&0&0&24.74&37.38&50.71\\ [.5ex]\hline
                \end{tabular}
		\label{table:ex1b}
	\end{table}
\begin{figure}[H]
\begin{center}
\includegraphics[width=120mm]{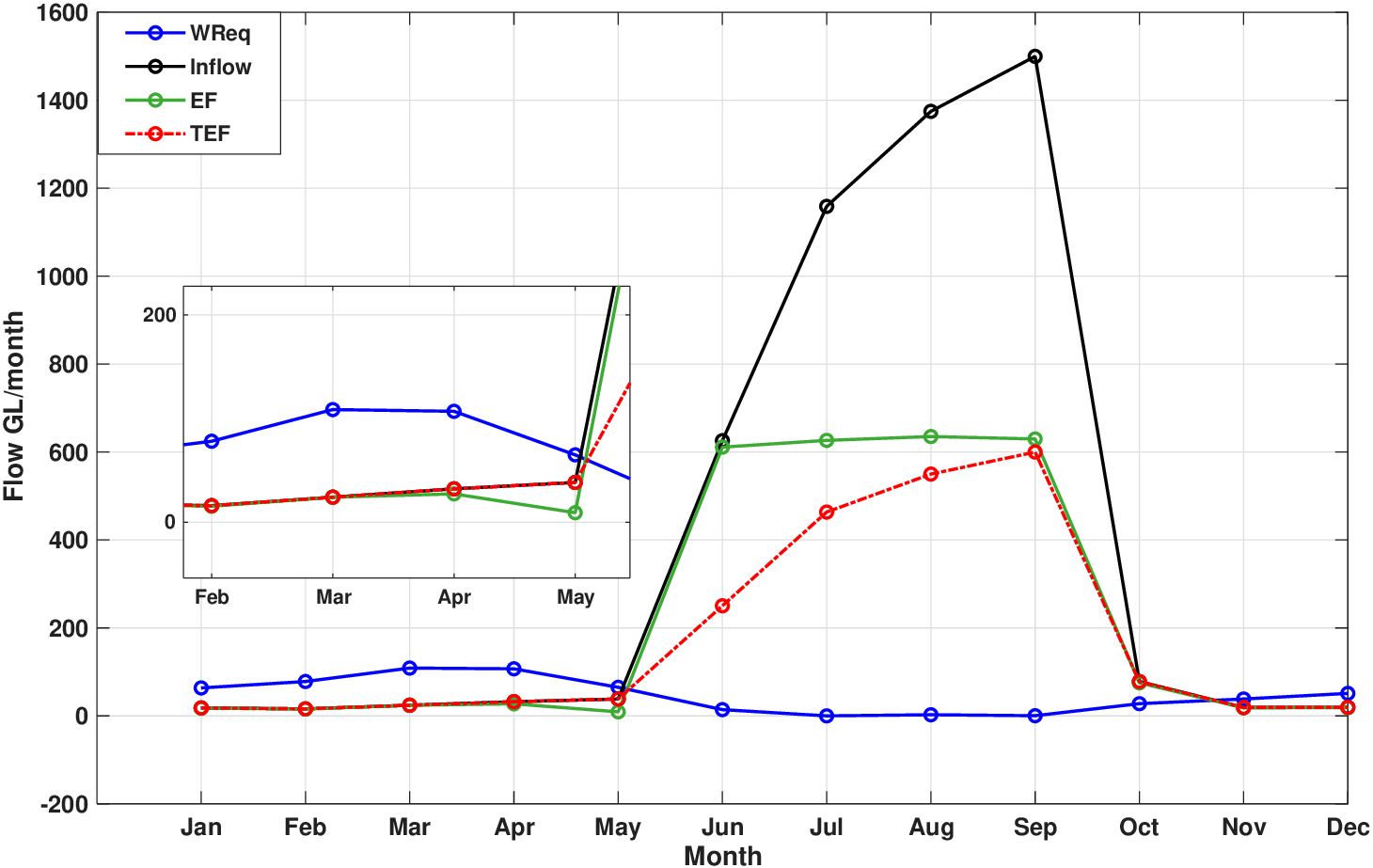} \\
\caption{Monthly Allocation of Environmental Flow and Water Requirement
for Model 2 in Dry Year.}
\label{fig_dry2}
\end{center}
\end{figure}
\noindent Table~\ref{table:ex1b} and Figure~\ref{fig_dry2}  display the optimal allocation for Model~2 in the dry season, where minimizing environmental flow deficiency is prioritized. Environmental flows are maintained throughout the year and increase during the monsoon months, while groundwater pumping supplements surface water availability, peaking at approximately 109~GL in March. As shown in Figure~\ref{fig_dry2}, the environmental flow (green) falls below the target environmental flow (red) during April, resulting in an environmental deficiency in this period. No deficiencies are observed in other months. Additionally, the optimal crop allocation satisfies only the minimum cultivation area requirements listed in Table \ref{table:mincrops}. This occurs because the model prioritizes minimizing environmental deficiency (reduces to $ 39.8460$), which leads to a reduction in net benefit to $ 1.7165 \times 10^{10}$ . 
\begin{table} [H]
 		\caption{\small{\textit{Optimal Allocation of Crops and Environmental Flow with Pumped Water for Model 1 in Average Year. Optimal Values are
        $f_1=-2.6785 \times 10^{10}$ and $f_2=451.6718$. }}}
		\footnotesize
		\vskip 1.5em
 		\centering
 		\begin{tabular}{|c| c| c| c| c| c| c| c| c| c| c| c|c|}
 			\hline
 			crops &   Aus &  Aman  & Boro & Wheat& Potato & Suger & Maize & Maize & Jute &  & &  \\ 
             &   rice &  rice  & rice & &  & cane & kharif 1& rabi &  &  &  &  \\ [0.5ex]
 			\hline
$X_c$&20000 &35000 & 30000 &10000&55271&16000&5000&5000&6000& & &\\ [.5ex]\hline
Env.&Jan &Feb & Mar &April&May&June&July&Aug&Sep&Oct&Nov & Dec\\ [.5ex]\hline
Flow&18 &0 & 0 &90&165 &3129&5793&6874&7500&385&65&51\\ [.5ex]\hline
$P_m$&0 &22.88 & 7.07 &0&0 &0&0&0&0&0&0&0\\ [.5ex]\hline
                \end{tabular}
		\label{table:ex2a}
	\end{table} 
\begin{figure}[H]
\begin{center}
\includegraphics[width=120mm]{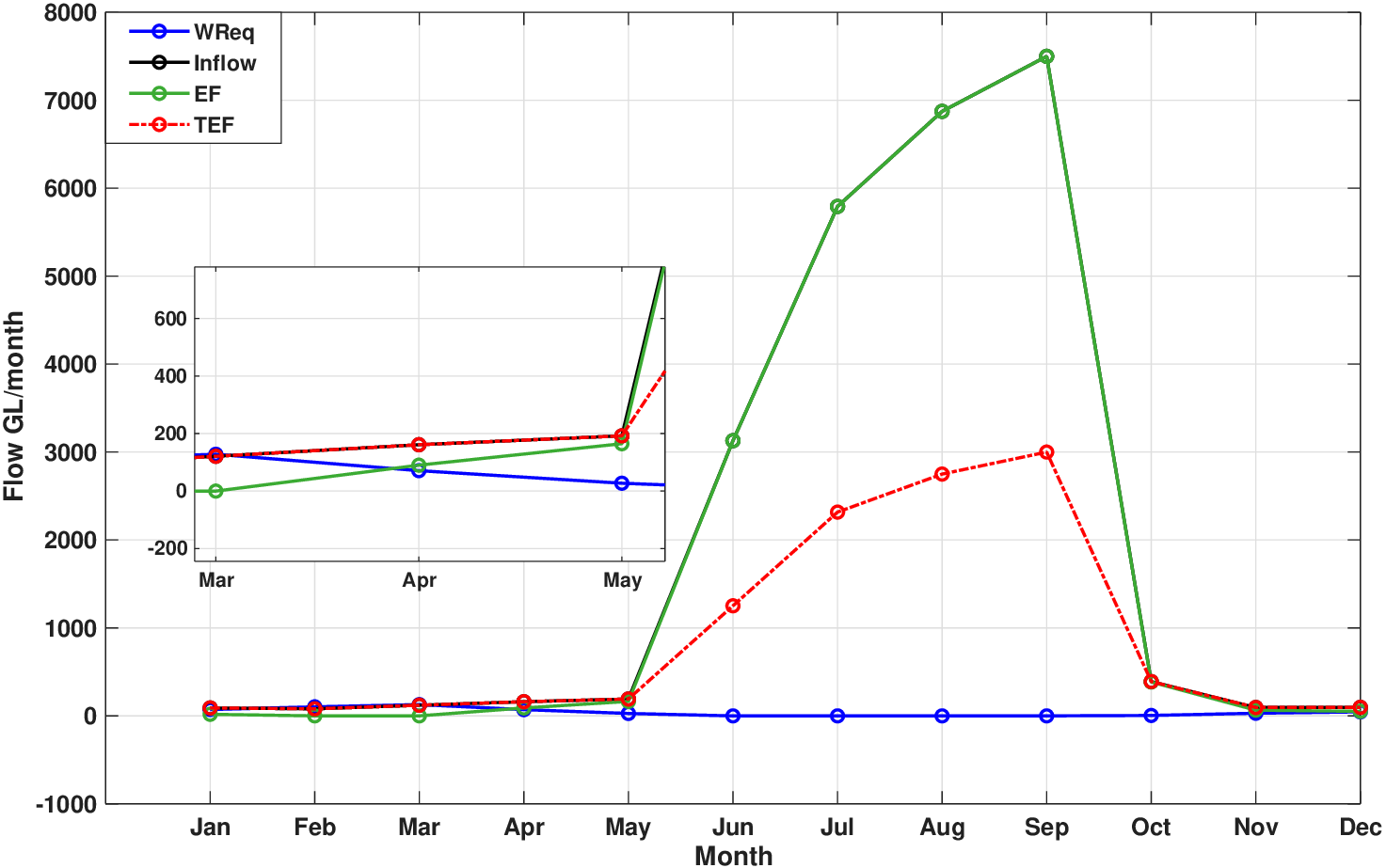} \\
\caption{{Monthly Allocation of Environmental Flow and Water Requirement for Model 1 in Average Year.}}
\label{fig_avg1}
\end{center}
\end{figure}
\noindent Table~\ref{table:ex2a} and Figure~\ref{fig_avg1} demonstrate that under net-benefit maximization in average hydrological years, potato receives the largest land allocation and the reason is described earlier in dry season scenerio. Environmental flow deficiency (EFD) occurs from November to May, reaching a maximum value of 451 GL during this period. Groundwater pumping is necessary only February and March. This is because the average year typically receives more rainfall than a dry year.

\begin{table} [H]
 		\caption{\small{\textit{Optimal Allocation of Crops and Environmental Flow with Pumped water for Model 2 in Average Year. Optimal Values are $f_1=-1.7098 \times 10^{10}$ and $f_2=0$.}}}
		\footnotesize
		\vskip 1.5em
 		\centering
 		\begin{tabular}{|c| c| c| c| c| c| c| c| c| c| c| c|c|}
 			\hline
 			crops &   Aus &  Aman  & Boro & Wheat& Potato & Suger & Maize & Maize & Jute &  & &  \\ 
             &   rice &  rice  & rice & &  & cane & kharif 1& rabi &  &  &  &  \\ [0.5ex]
 			\hline
$X_c$&20000 &35000 & 30000 &10000&15000&16000&5000&5000&6000& & &\\ [.5ex]\hline
Env.&Jan &Feb & Mar &April&May&June&July&Aug&Sep&Oct&Nov & Dec\\ [.5ex]\hline
Flow&90 &80 & 121 &161&192 &2827&2832&3838&6457&390&95&97\\ [.5ex]\hline
$P_m$&52.35 &63.38 & 87.88 &70.89&27.36 &0&0&0&0&4.76&29.93&41.4\\ [.5ex]\hline
                \end{tabular}
		\label{table:ex2b}
	\end{table}
\begin{figure}[H]
\begin{center}
\includegraphics[width=120mm]{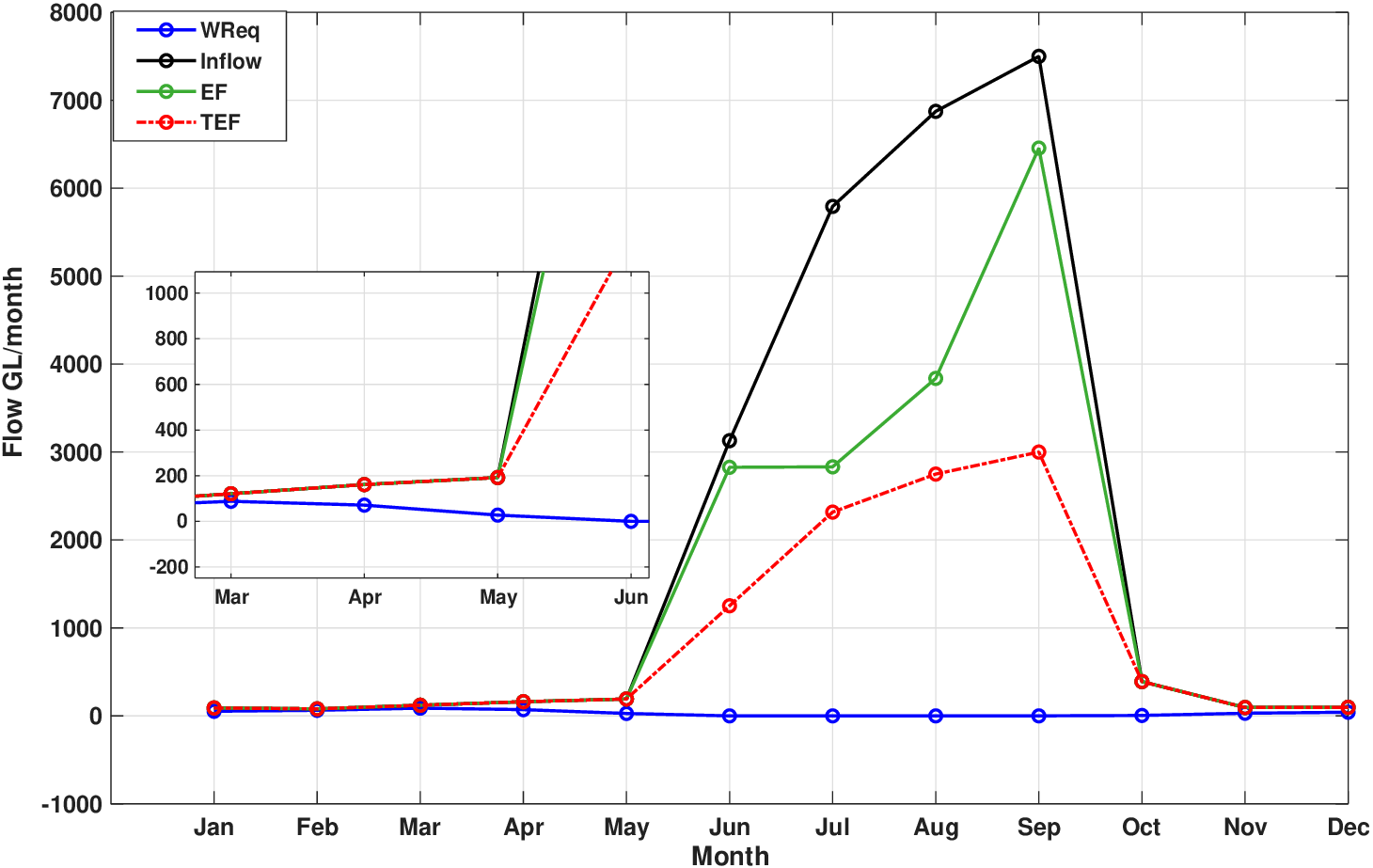} \\
 \caption{Monthly Allocation of Environmental Flow and Water Requirement
for Model 2 in Average Year.}
\label{fig_avg2}
\end{center}
\end{figure}

 \noindent Table~\ref{table:ex2b} and Figure~\ref{fig_avg2}  show that, under average hydrological conditions with the objective of minimizing environmental flow deficiency, environmental flows are prioritized, leading to minimum cropland allocation as the model restricts the use of the total cultivable land. There is no EFD. Environmental flows increase substantially during the monsoon season, and groundwater abstraction is used when necessary to meet these targets, leading to a reduction in net benefits.
\begin{table} [H]
 		\caption{\small{\textit{Optimal Allocation of Crops and Environmental Flow with Pumped Water for Model 1 in Wet Year. Optimal Values are $f_1=-2.6790 \times 10^{10}$ and $f_2=0.0000$. }}}
		\footnotesize
		\vskip 1.5em
 		\centering
 		\begin{tabular}{|c| c| c| c| c| c| c| c| c| c| c| c|c|}
 			\hline
 			crops &   Aus &  Aman  & Boro & Wheat& Potato & Suger & Maize & Maize & Jute &  & &  \\ 
             &   rice &  rice  & rice & &  & cane & kharif 1& rabi &  &  &  &  \\ [0.5ex]
 			\hline
$X_c$&20000 &35000 & 30000 &10000&55271&16000&5000&5000&6000& & &\\ [.5ex]\hline
Env.&Jan &Feb & Mar &April&May&June&July&Aug&Sep&Oct&Nov & Dec\\ [.5ex]\hline
Flow&162 &144 & 218 &290&346 &5634&10428&12374&13500&702&171&174\\ [.5ex]\hline
$P_m$&54.02 &73.8 & 84.06 &0&0 &0&0&0&0&0&6.06&29.19\\ [.5ex]\hline
                \end{tabular}
		\label{table:ex3a}
	\end{table}
\begin{figure}[H]
\begin{center}
\includegraphics[width=120mm]{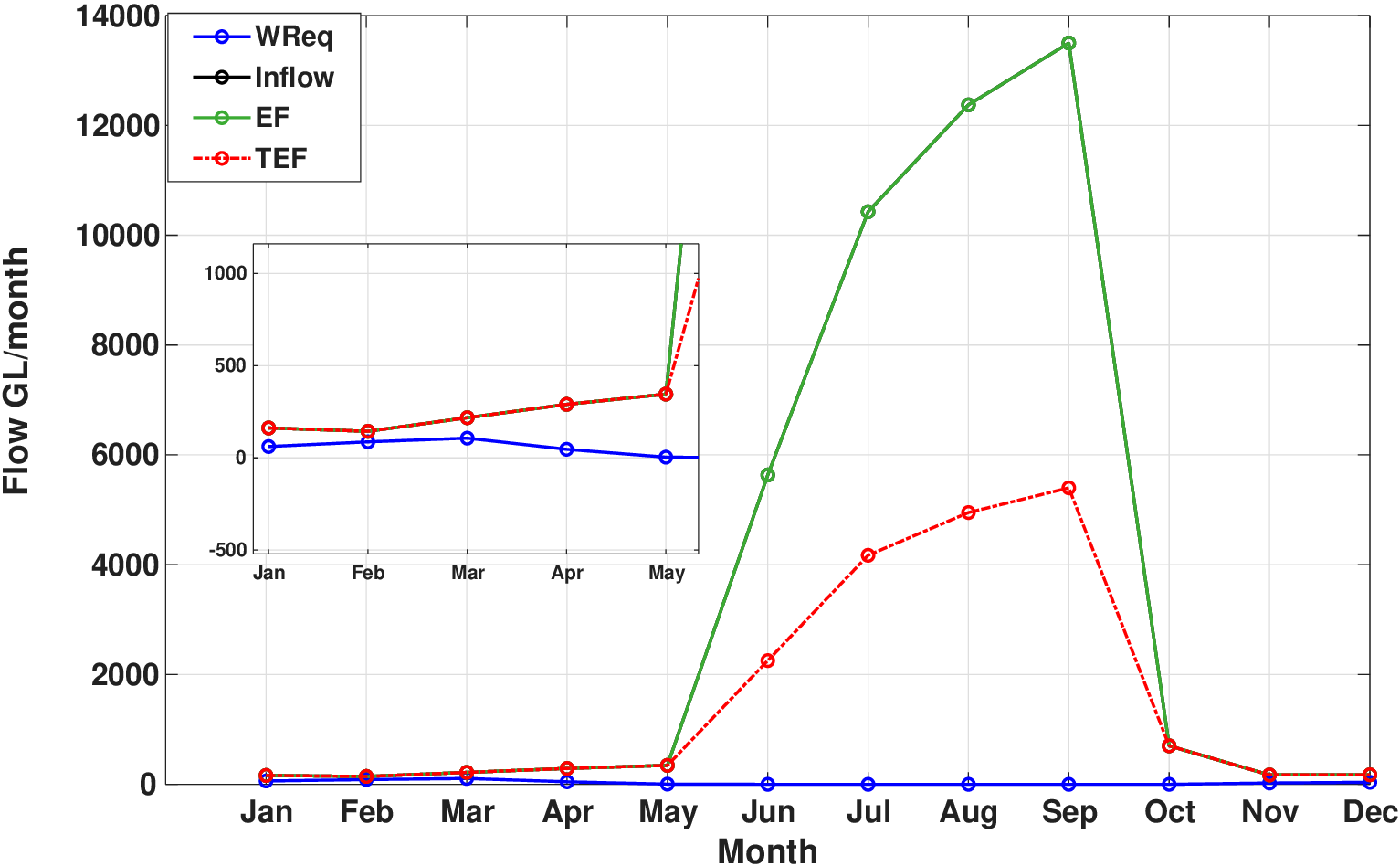} \\
\caption{{Monthly Allocation of Environmental Flow and Water Requirement for Model 1 in Wet Year.
}}
\label{fig_wet1}
\end{center}
\end{figure}
\noindent Table~\ref{table:ex3a} and Figure~\ref{fig_wet1} present that the same cropland allocation is observed under wet years, consistent with dry and average hydrological conditions, because land-use decisions are governed by crop productivity rather than water availability. Irrigation demands are met using combined surface and pumped water. There is no deficiency. Groundwater pumping is absent from April to October due to sufficient monsoon inflows, but increases from November to March under low rainfall and higher irrigation demand.

\begin{table} [H]
 		\caption{\small{\textit{Optimal Allocation of Crops and Environmental Flow with Pumped water for Model 2 in Wet Year. Optimal Values are $f_1=-1.6682 \times 10^{10}$ and $f_2=0$. }}}
		\footnotesize
		\vskip 1.5em
 		\centering
 		\begin{tabular}{|c| c| c| c| c| c| c| c| c| c| c| c|c|}
 			\hline
 			crops &   Aus &  Aman  & Boro & Wheat& Potato & Suger & Maize & Maize & Jute &  & &  \\ 
             &   rice &  rice  & rice & &  & cane & kharif 1& rabi &  &  &  &  \\ [0.5ex]
 			\hline
$X_c$&20000 &35000 & 30000 &10000&15000&16000&5000&5000&6000& & &\\ [.5ex]\hline
Env.&Jan &Feb & Mar &April&May&June&July&Aug&Sep&Oct&Nov & Dec\\ [.5ex]\hline
Flow&162 &144 & 218 &290&346 &2249&4428&6374&7500&702&171&174\\ [.5ex]\hline
$P_m$&44.95 &53.3 & 73.87 &46.29&3.71 &0&0&0&0&0&22.65&34.41\\ [.5ex]\hline
                \end{tabular}
		\label{table:ex3b}
	\end{table}
\begin{figure}[H]
\begin{center}
\includegraphics[width=120mm]{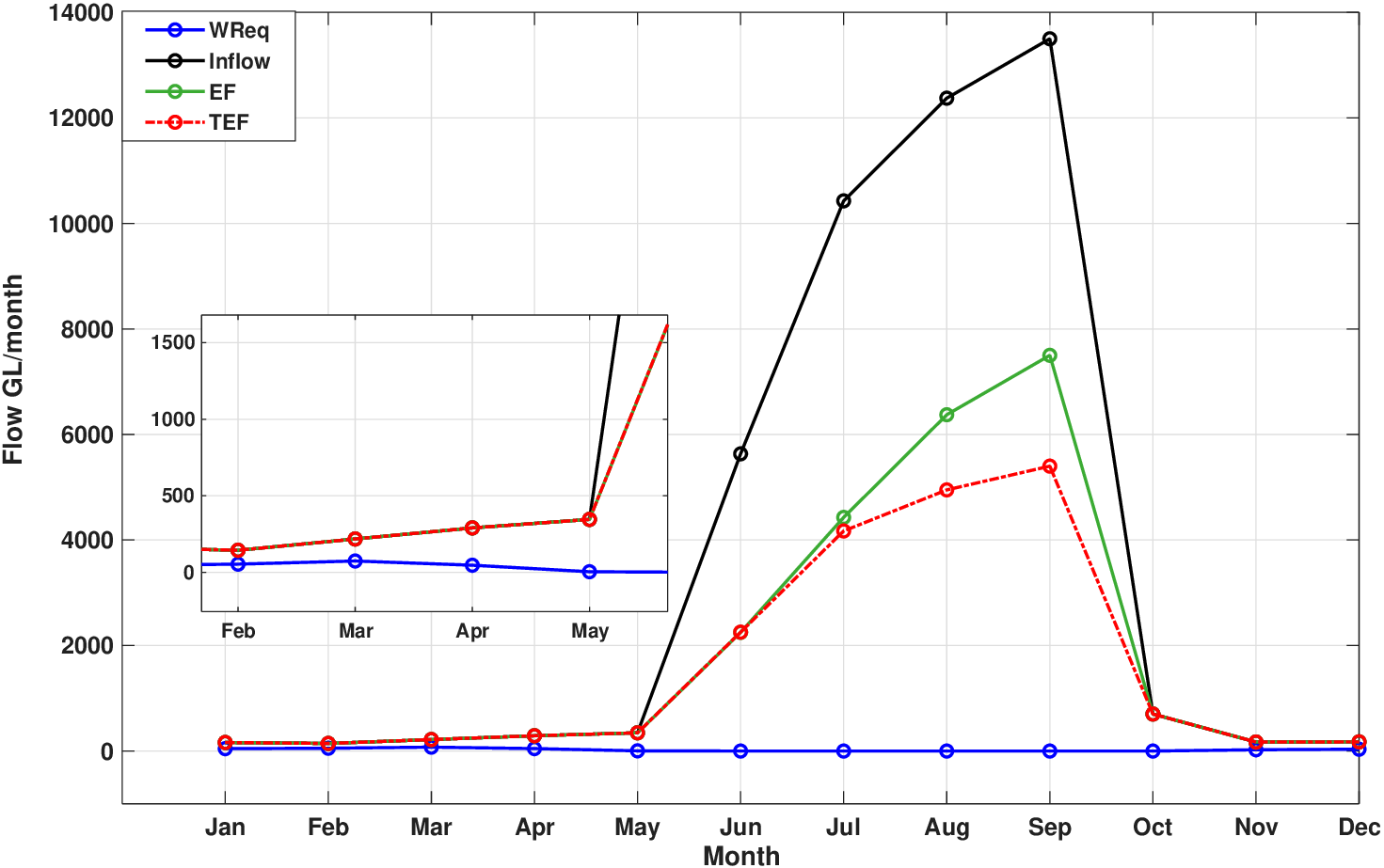} \\
 \caption{Monthly Allocation of Environmental Flow and Water Requirement
for Model 2 in Wet Year.}
\label{fig_wet2}
\end{center}
\end{figure}
\noindent Table~\ref{table:ex3b} and Figure~\ref{fig_wet2} depict that, when environmental sustainability is prioritized, cropland is allocated at a minimum level, limiting overall cultivation and reducing net benefit. This conservative land decision lowers irrigation demand and supports environmental goals. Consequently, surface water is reserved for environmental flows, groundwater is used in low-flow months, and monsoon rainfall eliminates pumping needs.
\par
\noindent Hence, the Models 1 and 2 described by \eqref{eqmathemodel1} and \eqref{eqmathemodel2}, respectively, are solved to determine the optimal net benefit and environmental flow deficit under alternative crop and environmental flow allocation in six different scenarios, as summarized in Table~\ref{table:ex1bc}.
\begin{table} [H]
 		\caption{\small{\textit{Optimal Environmental Flow Deficit months and Net Benefits Under Single-objective Optimization. }}}
		\footnotesize
		\vskip 1.5em
 		\centering
 		\begin{tabular}{|c| c| c| c| c| c| }
 			\hline
 			Years &   Single Obj. &  Month $EFD_{most}$  & Month $EFD_{least}$ & $NB_{max}$(Tk.) & $EFD_{min}$($GL$)\\ 
             [0.5ex]
 			\hline

$Dry$& max $f_1$ &May & June-Sep & $2.6746\times 10^{10}$ &194.9720\\ [.5ex]\hline
 &min $f_2$ &May &  June-March &$1.7165\times 10^{10}$&39.8460\\ [.5ex]\hline
Avg&max $f_1$ &March & May-Nov &$2.6785\times 10^{10}$&451.6718\\ [.5ex]\hline
 &min $f_2$ &- & Jan-Dec&$1.7098\times 10^{10}$&0\\ [.5ex]\hline
Wet&max $f_1$ & - & Jan-Dec &$2.6790\times 10^{10}$&0\\ [.5ex]\hline
 &min $f_2$ &- & Jan-Dec &$1.6682\times 10^{10}$&0\\ [.5ex]\hline
                \end{tabular}
		\label{table:ex1bc}
	\end{table}
 \par
\noindent The results show that the highest net benefit is reached  Tk. $2.6790 \times 10^{10}$ while the corresponding environmental flow deficit is reduced to zero, indicating no deficiency at the optimal solution.
\par
\noindent \textbf{Major Findings:}
The six scenarios yield different outcomes because each objective uniquely dictates water allocation between crops and the river within the given constraints:
\begin{itemize}

    \item \textbf{Groundwater pumping is structurally concentrated in late dry and early pre-monsoon months.}
    \par
\noindent Across all scenarios [Table~\ref{table:ex1a}-Table~\ref{table:ex3b}], February and March consistently experience insufficient surface water availability, making groundwater pumping the dominant source of irrigation during these months regardless of hydrological year type, due to low rainfall, high evapotranspiration, and peak irrigation water requirements.
    
\item \textbf{Net benefit maximization favors stable land allocation and flexible groundwater use.}  
\par
\noindent When the economic benefit is optimized, cropland allocation remains unchanged across dry, average, and wet years despite differences in inflow, rainfall, and evapotranspiration, because maintaining the most economically profitable cropping pattern by the crop productivity [Figure 2].
\item \textbf{Environmental flow deficiency minimization reduces net benefit.}  
\par
\noindent When environmental flow deficiency is prioritized (Model~2), a substantial share of surface inflow is reserved for ecological purposes [Figure~\ref{Inflow and TEF in Dry, Avg. and Wet years}], reducing the volume of water available for irrigation and increasing reliance on groundwater pumping to meet irrigation deficits. In addition, Model~2 allocates only the minimum required cropland to ensure environmental sustainability, which lowers total irrigated area and consequently results in reduced agricultural production and net economic returns compared to Model~1.

    \item \textbf{Optimized environmental flows drive monthly groundwater pumping differences.}
    \par
\noindent Monthly  water requirement is determined by the optimized allocation of crop areas. Whether this requirement is met using surface water or groundwater depends on the optimized environmental flow releases. As a result, even within the same hydrological year, Model 1 and Model 2 show different month-by-month groundwater pumping patterns and total pumping volumes [Table~\ref{table:ex1a}-Table~\ref{table:ex3b}]. These differences arise because monthly inflow, rainfall, and evapotranspiration interact differently [Figure~\ref{Inflow and TEF in Dry, Avg. and Wet years} - Figure~\ref{fig_R_m & ET_m}]
    
\item \textbf{Peak inflows sustain irrigation without groundwater pumping, favoring short-duration crops.} 
\par
\noindent Across all scenarios Figure~[\ref{fig_dry1}]-Figure~[\ref{fig_wet2}] peak inflows occur during  May to October. During this period, surface water is sufficient to meet irrigation demands, eliminating the need for groundwater pumping. Farmers in Rajshahi, therefore, cultivate very short-duration crops, such as vegetables, mung beans, Black gram, and sesame [Section \ref{Case Study in Rajshahi}].
\end{itemize}



\subsection{Multiobjective Optimization Problems} \label{MOP}

\noindent In this section, we aim to optimise two competing objectives: maximising net benefit while simultaneously minimising the environmental impact of water distribution. To achieve this, we develop a multiobjective model. Multiobjective optimisation is designed to handle trade-offs between goals that often conflict with one another. Since we cannot maximise profit without impacting the environment, our goal is to approximate a set of "trade-off" solutions, known as Pareto or weak Pareto solutions. These represent non-dominated outcomes where one objective cannot be improved without worsening the other. Essentially, we are looking for the optimal balance between economic profit and environmental sustainability.

\noindent To solve this, we employ scalarization methods. These techniques transform the complex multiobjective problem into a standard single-objective problem, allowing us to use traditional optimisation algorithms. This process involves using normalisation parameters to combine the two objective functions. Parameter defined by the relation below.

\begin{equation}
    \hspace*{20mm} W := \left
\{ w \in \mathbb{R}^{2} \mid w_1, w_2 > 0,\sum_{i=1}^{2} w_i=1
\right \}.
\end{equation}

\noindent We utilise the weighted-constraint method introduced by Burachik et al. \cite{Burachik2014}. We prefer this approach for its ease of implementation. For readers interested in further details on this and other scalarization approaches, refer to \cite{Burachik2014}, \cite{BurKayRiz2017}, \cite{BurKayRiz2022}.

\noindent In the weighted-constraint method, two subproblems are constructed since Model $3$ is a bi-objective model. In the first subproblem, the net benefit of crop production is treated as the objective function, while the second objective is incorporated into the constraint set \(X\). Given that MATLAB is used for minimization, the first subproblem is formulated by minimizing \(-f_1\) subject to the extended constraint \(w_2 f_2 \leq w_1 f_1\), together with the original constraints in \(X\). A similar procedure is applied in the second subproblem. In this case, minimizing environmental deficiency \((f_2)\) is taken as the objective function, while \(f_1\) is included in the constraint set.

\noindent Let $w \in W$, \textbf{\em the weighted-constraint approach} is defined as follows.

\noindent 
\begin{equation} \label{nchep1}
\mbox{Subproblem 1:}\ \left\{\begin{array}{rl} \ds\min_{(X_c, Env.Flow_m) \in  X} & \
\ - w_1f_1(X_c, Env.Flow_m),
\\[4mm]
\mbox{subject to} & \ \ w_2f_2(X_c, Env.Flow_m) \leq w_1f_1(X_c, Env.Flow_m),\\ & \ \ \ds X_c, Env.Flow_m  \geq 0.
\end{array}
\right.
\end{equation}
and
\begin{equation} \label{nchep2}
\mbox{Subproblem 2:}\ \left\{\begin{array}{rl} \ds\min_{(X_c, Env.Flow_m)\in  X} & \
\ w_{2}f_{2}(X_c, Env.Flow_m),
\\[4mm]
\mbox{subject to} & \ \ w_1f_1(X_c, Env.Flow_m) \leq w_2f_2(X_c, Env.Flow_m), \\ & \ \ \ds X_c, Env.Flow_m  \geq 0.
\end{array}
\right.
\end{equation}
where $X$ is defined in \eqref{eqmathemodel1}, $f_1$ in \eqref{revobj1} and $f_2$ in \eqref{envobj}.

\noindent In Case Study \ref{exp1}, net profit and environmental deficiency are optimized separately, and the results identify the optimal solutions when each objective is optimized individually subject to the given constraints. We now turn to the trade-off between these two objectives. In this situation, it is difficult to find a single solution that optimizes both objectives simultaneously. Instead, we approximate a set of non-dominated solutions that represent the trade-off between the two objectives; this set is known as the Pareto front. In Case Study \ref{exp2}, we solve the multiobjective optimization problem defined in \eqref{eqmathemodel3}.
\color{black}
\begin{casestudy}\label{exp2} 
In this study, we address the multiobjective optimization model defined in \eqref{eqmathemodel3}, using the same study area, parameters, and constraints for the Rajshahi zone discussed in Case Study \ref{exp1}. 
\end{casestudy}
\noindent All subproblems described in \eqref{nchep1} and \eqref{nchep2} are solved using the fmincon solver with range of algorithms, as in Case Study~\ref{exp1}. The resulting Pareto fronts for the three hydrological periods are shown in  Figure~\ref{fig_wet2}. The algorithmic steps are detailed in the Appendix \ref{mopalg}.

\begin{figure}[H]
\begin{center}
\includegraphics[width=120mm]{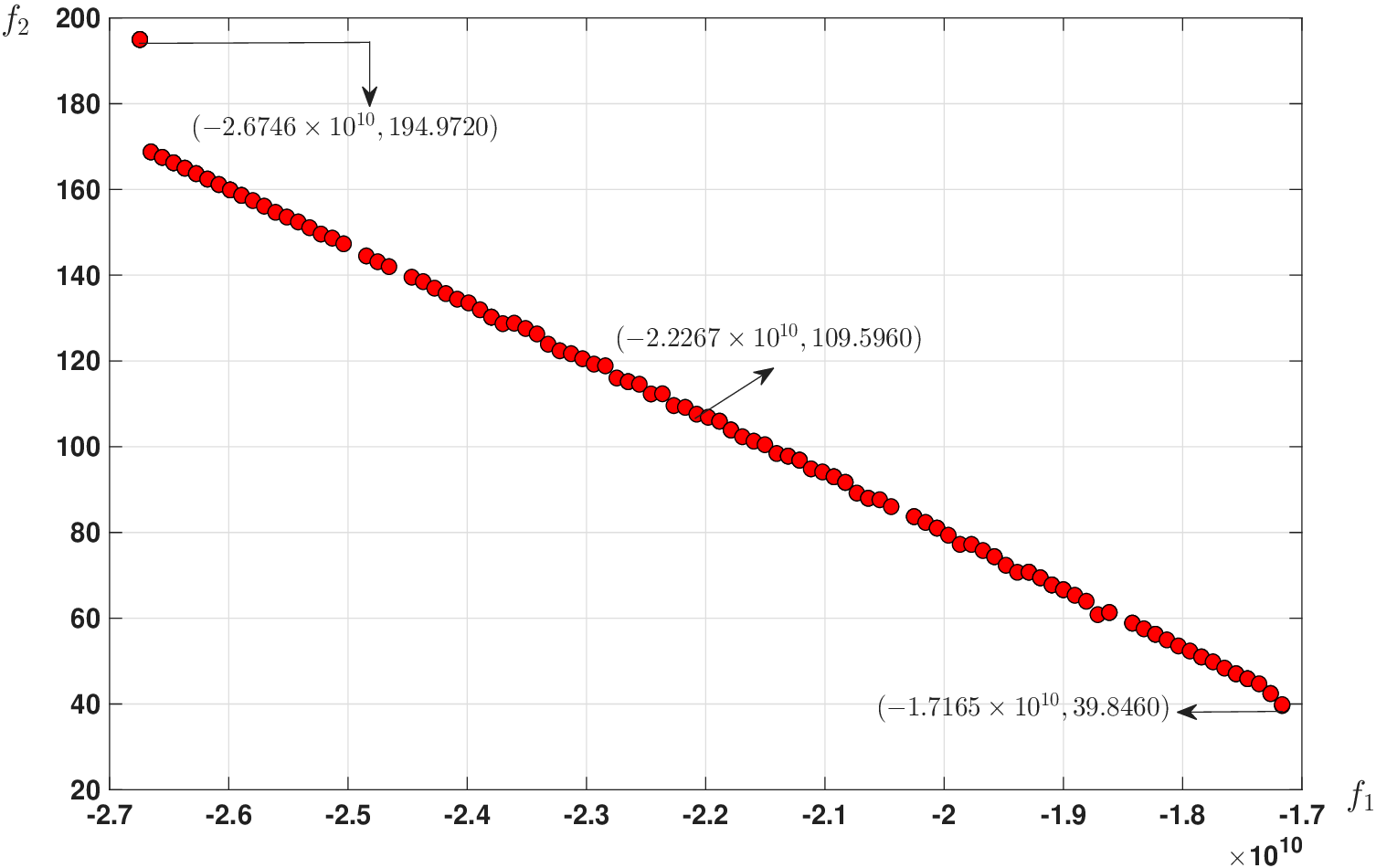} \\
\caption{Pareto front of Multiobjective Approach in Dry Year}
\label{mopfig1}
\end{center}
\end{figure}

\begin{figure}[H]
\begin{center}
\includegraphics[width=120mm]{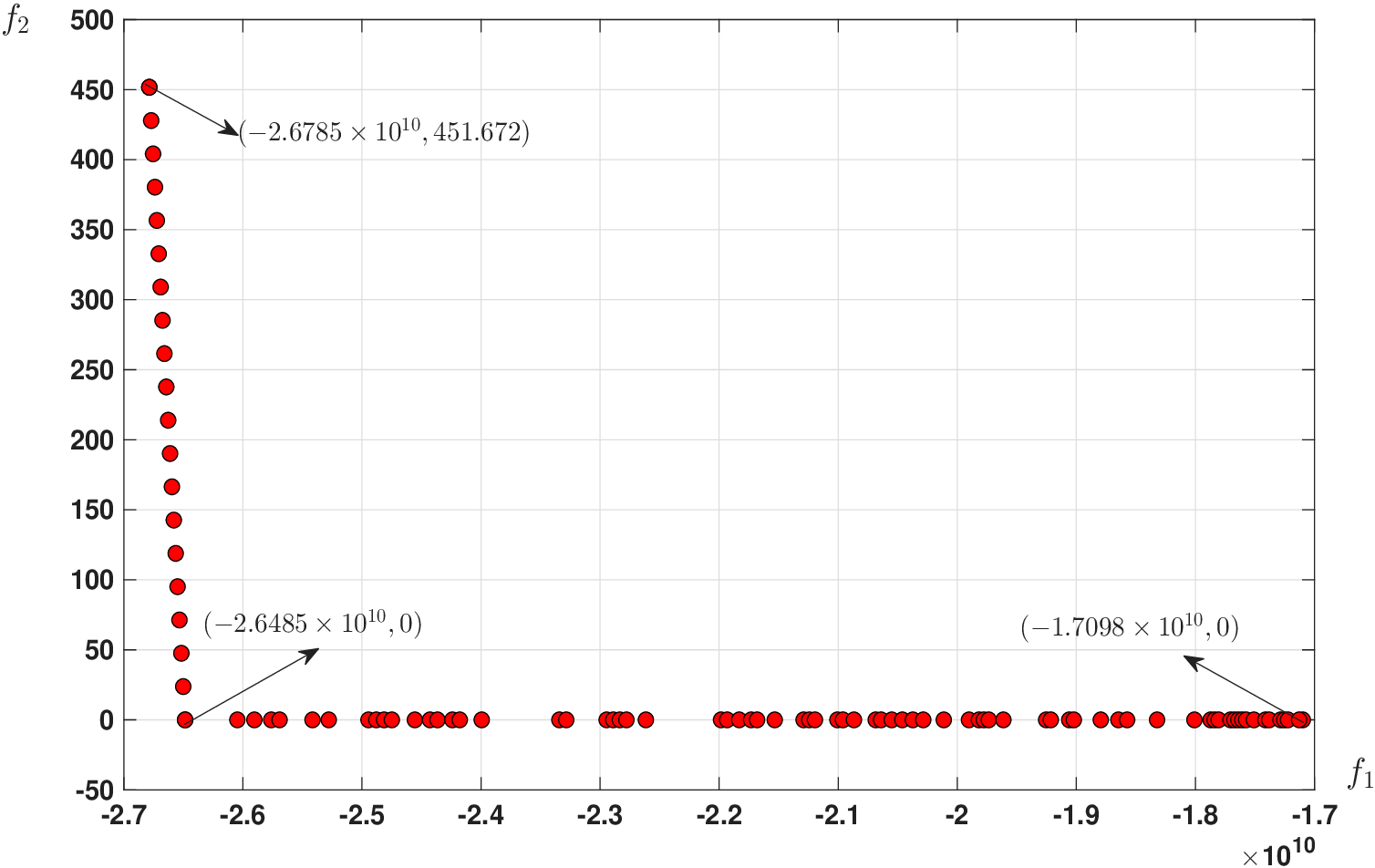} \\
\caption{Pareto front of Multiobjective Approach in Average Year}
\label{mopfig2}
\end{center}
\end{figure}
\begin{figure}[H]
\begin{center}
\includegraphics[width=120mm]{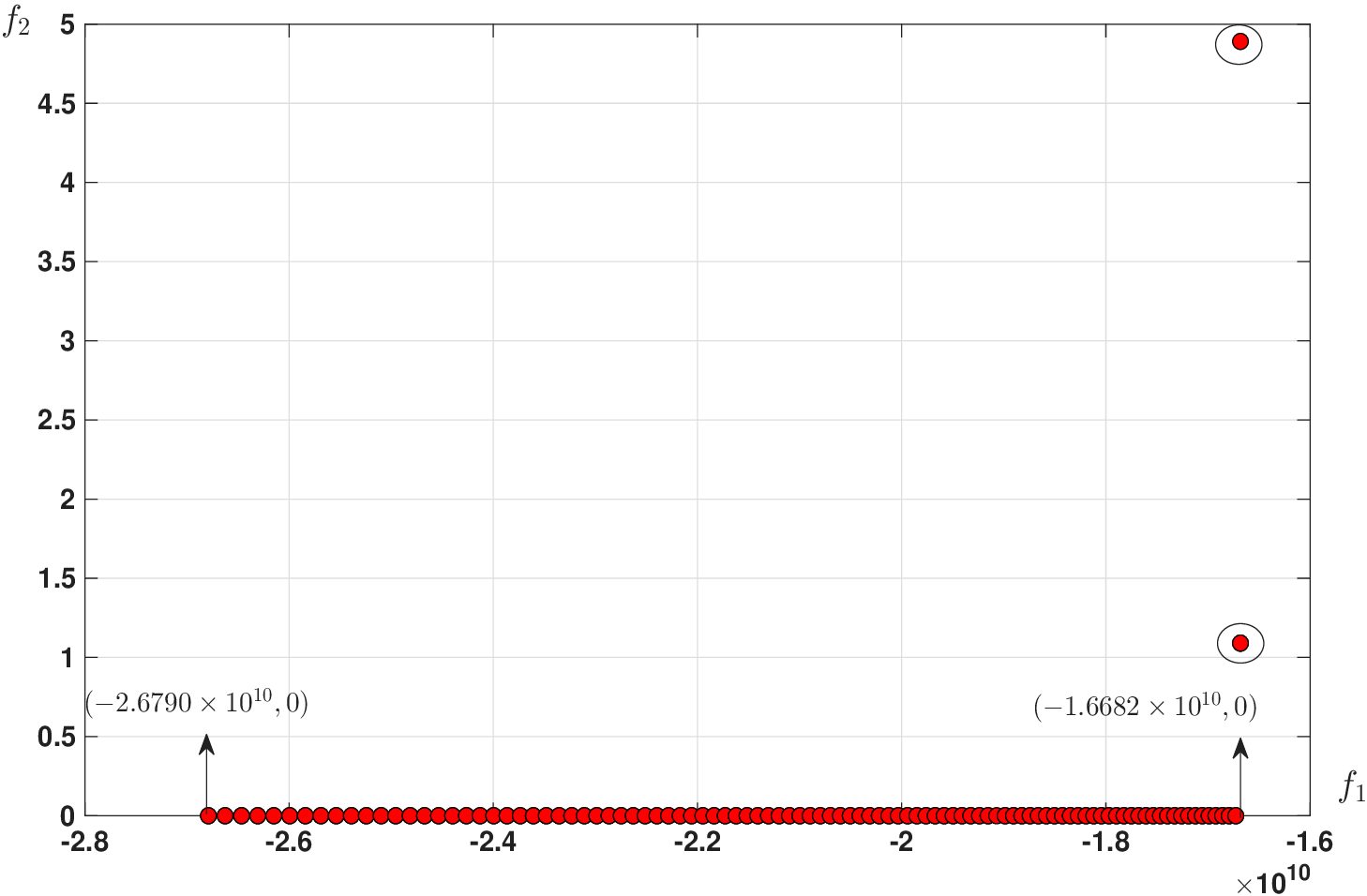} \\
\caption{Pareto front of Multiobjective Approach in Wet Year}
\label{mopfig3}
\end{center}
\end{figure}
\par
\noindent \textbf{Comparison Analysis:} In this analysis, we solve the multiobjective optimization problem using the `gamultiobj` and `paretosearch` algorithms available in MATLAB, and we additionally apply a scalarization-based weighted-constraint approach as formulated in equations \eqref{nchep1} and \eqref{nchep2}. The performance of Model~3 is evaluated and compared across all methods under identical parameter settings.

\noindent The results indicate that `paretosearch` is largely unable to generate a representative Pareto front for this problem. While `gamultiobj` performs slightly better, it is only capable of approximating a portion of the Pareto front. In contrast, the weighted-constraint method successfully produces Pareto-optimal solutions, provided that a sufficiently good initial guess is supplied when solving equations \eqref{nchep1} and \eqref{nchep2}. This happens because fmincon uses gradient-based methods that often converge to a local optimum. To reduce this issue, we randomize the initial starting points and run the model multiple times, which helps achieve better coverage and more reliable approximations of the Pareto front.

\noindent Figure~\ref{mopfig1} illustrates the Pareto front for the dry year, representing a scenario in which economic performance is dominant. The results show that maximizing net benefit leads to a substantial compromise in environmental flow, with the deficiency increasing to approximately $194.97$ GL, while the net benefit reaches Tk. $(2.6746 \times 10^{10})$. In contrast, giving priority to the environmental objective significantly reduces the flow deficiency to \(39.8460\) GL, but at the cost of a lower net benefit of Tk. $(1.1765 \times 10^{10})$. This trade-off clearly demonstrates the inherent conflict between economic returns and ecological sustainability.

\noindent  Figure~\ref{mopfig2} shows the Pareto front for the average year. Here, net benefit ranges from Tk.$1.7098\times 10^{10}$ to Tk. $2.6486\times 10^{10}$ without impacting environmental flow requirements; however, the environmental deficiency increases gradually and reaches approximately 452 GL at the upper bound of net benefit Tk.$2.6785\times 10^{10}$. This outcome demonstrates that moderate profit gains may be achieved without severe ecological penalties, but aggressive profit maximization imposes substantial environmental costs. Such results provide actionable insight for farmers and local water authorities seeking balanced management strategies.

\noindent  Figure~\ref{mopfig3} reflects a hydrological year with higher rainfall and lower evapotranspiration. Under these conditions, the Pareto front becomes comparatively weak and no environmental deficiency is observed over a profit interval of Tk. $1.6682\times 10^{10}$ to Tk. $2.6790\times 10^{10}$. The two identified outlier points are due to solver difficulties, which are marked accordingly. This scenario underscores the influence of climatic conditions on the shape and severity of the trade-off, reinforcing that environmental constraints may relax substantially in wetter years.

\section{Conclusion}\label{Con}

\noindent We propose a multiobjective mathematical modeling framework for agricultural water allocation that integrates canal capacity, crop-specific minimum area requirements, and non-negative allocation conditions for both surface and groundwater use. The model aims to achieve an optimal balance between economic efficiency and environmental sustainability across varying hydrological years—dry, average, and wet. By incorporating realistic constraints such as canal flow limits and minimum irrigation demands per crop, the framework overcomes common simplifications in previous studies, which often neglected operational feasibility or assumed homogeneous land use patterns.

\noindent  To validate the proposed framework, a series of computational experiments were conducted under distinct inflow–environmental flow scenarios, solved using Sequential Quadratic Programming (SQP) and compared with existing metaheuristic approaches. The results demonstrate that the SQP-based scalarization approach significantly reduces computation time while maintaining high accuracy in Pareto front generation. This provides a reliable alternative to population-based methods such as NSGA-II, especially for decision-making contexts where quick, interpretable results are necessary.

\noindent The developed model not only determines the optimal allocation of surface and groundwater among crops but also evaluates trade-offs between total benefit maximization and environmental flow satisfaction. Sensitivity analyses reveal how changes in canal capacity, inflow variability, and environmental flow thresholds affect both economic returns and ecological outcomes. In particular, scenarios involving reduced inflow conditions show that the model effectively reallocates groundwater while maintaining the minimum required environmental flow, ensuring sustainability under scarcity.

\noindent While the proposed model demonstrates strong applicability, certain limitations remain. The current framework assumes deterministic hydrological and agronomic parameters, which may not fully capture stochastic inflow fluctuations or demand uncertainties in real-world systems. Moreover, energy and pumping cost variations have not been explicitly modeled. Future research could extend this work toward stochastic multiobjective optimization by integrating uncertainty in inflow, crop demand, and groundwater recharge, as well as environmental factors such as carbon footprint and energy consumption. Incorporating dynamic canal operation rules and hybrid water sources (e.g., reclaimed or rainwater) would further enhance its adaptability to changing climatic and agricultural conditions.

\noindent \textbf{\large {Acknowledgment:}} 
\noindent We acknowledge Chatgpt(5.2) AI tools, which are used to rephrase, edit and polish the author's written text for spelling, grammar, or general style. 

\section{\small Funding and/or Conflicts of interests/Competing interests}
 The authors declare that no funds, grants, or other support were received during the preparation of this manuscript. The authors have no relevant financial or non-financial interests to disclose.



\section{Appendix}\label{appen}

\subsection{Algorithm Steps for Models 1 and 2}\label{solprocedure}
\noindent The solution procedure for Models $1$ and $2$ is summarized below, consistent with the approach described in Case Study \ref{exp1}.
 \begin{description}
\item[Step $\mathbf{1}$] { \textbf{(Input)}} \\ 
In this step, all relevant parameters and input data for the optimization model are retrieved from the dataset and defined. A detailed description of the parameters is provided in the Case Study \ref{exp1} definition.
 \item[Step  $\mathbf{2}$] {\textbf{(Outline the Problem)}}\\
The constrained optimization problems are formulated as in \eqref{eqmathemodel1} for maximizing revenue and in \eqref{eqmathemodel2} for minimizing environmental flow deficiency. These problems are solved using MATLAB’s `fmincon` solver with default settings and randomized initial points. Four algorithm options are employed: Active Set, Interior-Point, SQP, and SQP-Legacy. The lower bounds are set based on the available data, while the upper bounds are left unconstrained for both $X_c$ and $Env.Flow_m$.

 \item[Step  $\mathbf{3}$] \textbf{(Solve Problem)}\\
Solve $\ds \max f_1$ and $\ds \min f_2$ under the constraints set $X$, and record the solutions $\left(\bar{X^1}_{c}, \bar{Env.Flow^1}_{m}\right)$.

 \item [Step  $\mathbf{4}$] \textbf{(Algorithm Validation)}\\
Multiple algorithms are implemented using fmincon, and their performance is evaluated. The results are then compared and recorded.
\end{description}

\subsection{Algorithm Steps for Model~3}\label{mopalg}
\noindent The following algorithm is used in Case Study-\ref{exp2} to generate the Pareto fronts for objectives $f_1$ and $f_2$, respectively.

\noindent \textbf{Algorithm 1} \label{algo1}
\begin{description}
	\item[Step $\mathbf{1}$] { \textbf{(Input)}} \\ Assign all parameters as specified in Case Study \ref{exp1}.
	\item[Step  $\mathbf{2}$] 	
\textbf{(Determine the individual optima)} Solve the problems \(\min f_1\) and \(\min f_2\) subject to the constraints defined in \eqref{eqmathemodel3}, using the parameter values specified in Case Study \ref{exp1}. The corresponding solutions are denoted by \(\left(\bar{X}_{c}, \bar{Env.Flow}_{m}\right)\).

	\item[Step $\mathbf{3}$] { \textbf{(Generate weighted parameters)}} \\ We generate weight vectors $w_i$ in this step. We use weight generation method introduced in \cite[Step 3 of Algorithm 3]{BurKayRiz2022}.
	\item[Step $\mathbf{4}$] From Step`3, we obtain $w=(w_1,1-w_1)$.
	\begin{description}
		\item[(a)] Solve Subproblems \eqref{nchep1} and  \eqref{nchep2}. Record Solutions $\bar{x}_1:=\left(\bar{X^1}_{c}, \bar{Env.Flow^1}_{m}\right)$ and $\bar{x}_2:=\left(\bar{X^2}_{c}, \bar{Env.Flow^2}_{m}\right)$, respectively. 
		\item[(b)] Find non-dominated points : 
				\begin{description}
			\item[(i)] If $\bar{x}_1=\bar{x}_2$, then record the point  $\bar{x}=\bar{x}_1$ as a non-dominated point.
			\item [(ii)] If $\bar{x}_1\neq \bar{x}_2$, then discard the dominated point. \\
			Record the data.	
		\end{description}
	\end{description}
	\item[Step $\mathbf{5}$] (Output)\\
	All recorded points are Pareto points of Case Study-\ref{exp2}.\\
\end{description}

\subsection{Data}\label{data}
\par
\noindent The data on prices, yields, and variable costs for selected crops in Rajshahi were collected from agriculture related online journal portals, including ais.rajshahi.gov.bd, bssnews.net, observerbd.com, shomoyeralo.com, barciknews.com, and sarabangla.net.  and crop coefficients are collected from the sources \cite{ Hossain2020, Mila2016, Milla2018, Mojid2020, Rahman2015, Sen2019}, and the weather data is collected from \cite{ BWDB2021,HAli2021}.
\begin{table} [H]
 		\caption{\small{\textit{Price (Tk/ton), Yield (ton/ha), and Variable Cost (Tk. /ha) for Selected Crops}}}
		\footnotesize
		\vskip 1.5em
 		\centering
 		\begin{tabular}{|c| c| c| c| c| c| c| c| c| c|}
 			\hline
 			crops &   Aus &  Aman  & Boro & Wheat& Potato & Suger & Maize & Maize & Jute\\ 
             &   rice &  rice  & rice & &  & cane & Kharif 1& Rabi &\\ [0.5ex]
 			\hline
$P_c$&33000 &28000 & 36000 &34000&13000&5580&32000&3100&8600\\ [.5ex]\hline
$Y_c$&5.96	&4.44	&6.2	&3.04	&20.44	&60	&4.5	&5.5	&3.05\\ [.5ex]\hline
$Vcost_c$&80610	&52500	&132900	&46995	&28468	&132900	&43260	&73500	&35790\\ [.5ex]\hline

                \end{tabular}
		\label{table:ex50}
	\end{table}
    \begin{table} [H]
 		\caption{\small{\textit{Crop Coefficient ($K_{c,m})$ {}}}}
		\footnotesize
		\vskip 1.5em
 		\centering
 		\begin{tabular}{|c| c| c| c| c| c| c| c| c| c| c| c|c|}
 			\hline
 			Month.&Jan &Feb & Mar &April&May&June&July&Aug&Sep&Oct&Nov & Dec\\ [.5ex]\hline 
$Aus Rice$ & 0 & 0 & 0 & 1.05 & 1.10 & 1.20 & 1.20 & 0.90 & 0 & 0 & 0 & 0 \\[.5ex]\hline
$Aman Rice$   & 0 & 0 & 0 & 0 & 0 & 1.05 & 1.07 & 1.07 & 1.20 & 0.87 & 0 & 0 \\[.5ex]\hline
$Boro Rice$   & 1.10 & 1.16 & 1.16 & 0.86 & 0 & 0 & 0 & 0 & 0 & 0 & 0 & 1.05 \\[.5ex]\hline
$Wheat$       & 1.12 & 0.48 & 0 & 0 & 0 & 0 & 0 & 0 & 0 & 0.42 & 0.78 & 0 \\[.5ex]\hline
$Potato$      & 0.62 & 1.15 & 0.80 & 0 & 0 & 0 & 0 & 0 & 0 & 0 & 0 & 0.25 \\[.5ex]\hline
$Sugarcane$   & 0.40 & 1.15 & 1.15 & 1.25 & 1.25 & 1.25 & 1.25 & 1.25 & 1.25 & 0.75 & 0.75 & 0.75 \\[.5ex]\hline
$Maize_{Kharif-1}$& 0 & 0 & 0 & 0.64 & 1.13 & 1.13 & 0.66 & 0 & 0 & 0 & 0 & 0 \\[.5ex]\hline
$Maize_{Rabi}$ & 0.75 & 0 & 0 & 0 & 0 & 0 & 0 & 0 & 0 & 0.38 & 0.87 & 1.36 \\[.5ex]\hline
$Jute$        & 0 & 0 & 0 & 0.72 & 1.39 & 1.26 & 0.46 & 0 & 0 & 0 & 0 & 0 \\[.5ex]\hline

             \end{tabular}
		\label{table:ex5b}
	\end{table}

\begin{table} [H]
 		\caption{\small{\textit{Rainfall and Evapotranspiration  ($ 10^{-4}$GL/ha)}}}
		\footnotesize
		\vskip 1.5em
 		\centering
 		\begin{tabular}{|c| c| c| c| c| c| c| c| c| c| c| c|c|}
 			\hline
 			Month.&Jan &Feb & Mar &April&May&June&July&Aug&Sep&Oct&Nov & Dec\\ [.5ex]\hline 
            
$Rain_{Dry}$&0.6	&0.97	&1.6	&4.44	&9.64	&15.65	&20.55	&16.31	&17	&7.69	&0.83	&0.60
\\ [.5ex]\hline
$Rain_{Avg}$&0.9	&1.46	&2.42	&6.7	&14.55	&23.62	&31.02	&24.62	&25.66	&11.6	&1.263	&0.9
\\ [.5ex]\hline
$Rain_{Wet}$&1.11	&1.82	&3.02	&8.36	&18.15	&29.46	&38.69	&30.71	&32	&14.47	&15.7	&11.2
\\ [.5ex]\hline
$ET_{Dry}$&10.7	&11.3	&17.8	&19	&19.7	&15.5	&14.3	&14.3	&13.8	&12.4	&8.6	&9
\\ [.5ex]\hline
$ET_{Avg}$&9.3	&9.8	&15.5	&16.5	&17.1	&13.5	&12.4	&12.4	&12	&10.8	&7.5	&7.8
\\ [.5ex]\hline
$ET_{Wet}$&8.4	&8.8	&14	&1.49	&15.4	&12.2	&11.2	&11.2	&10.8	&9.7	&6.8	&7
\\ [.5ex]\hline

                \end{tabular}
		\label{table:ex5c}
	\end{table}

\noindent \textbf{Author Contributions}: 
All authors contributed equally to this study.
  
\end{document}